\input ./style/arxiv-general.cfg
\documentclass[aap,MSNbibl,seceqn,dvips]{arximspdf}
\makeatletter
   \@ifpackageloaded{graphicx}{}{\usepackage{graphicx}}
\makeatother

\doi{10.1214/15-AAP1099}
\volume{26}
\issue{1}
\pubyear{2016}
\firstpage{507}
\lastpage{548}
\docsubty{FLA}

\makeatletter
\newcommand{\rrvert}{\vert}
\newcommand{\llvert}{\vert}
\newcommand{\eqref}[1]{(\ref{#1})}
\newtheorem{thmm}{Theorem}[section]
\newtheorem{lem}[thmm]{Lemma}
\newtheorem{lemm}{Lemma}[section]

\newtheorem{prop}[thmm]{Proposition}

\newproclaim{rem}[thmm]{Remark}
\newproclaim{defn}[thmm]{Definition}
\newproclaim{ex}[thmm]{Example}
\newproclaim{conjec}[thmm]{Conjecture}

\newcommand{\I}{\mathrm{i}}
\newcommand{\pp}{\mathbb{P}}
\newcommand{\ee}{\mathbb{E}}
\newcommand{\rr}{\mathbb{R}}
\newcommand{\nn}{\mathbb{N}}
\newcommand{\zz}{\mathbb{Z}}
\newcommand{\cc}{{\mathbb{C}}}
\newcommand{\ep}{\varepsilon}
\newcommand{\wt}{\widetilde}
\newcommand{\mfh}{\mathfrak{h}}
\newcommand{\mff}{\mathfrak{f}}
\newcommand{\mfg}{\mathfrak{g}}

\newcommand{\Res}{\operatorname{Res}}

\let\nnu\nu

\newcommand{\ts}{\widetilde{s}}
\newcommand{\tw}{\widetilde{w}}

\makeatother

\begin{document}
\begin{frontmatter}

\title{Exact formulas for random growth with half-flat initial data}
\runtitle{Exact formulas for half-flat random growth}

\begin{aug}
\author[A]{\fnms{Janosch}~\snm{Ortmann}\ead[label=e1]{janosch.ortmann@utoronto.ca}\thanksref{T1,m1}},
\author[A]{\fnms{Jeremy}~\snm{Quastel}\ead[label=e2]{quastel@math.toronto.edu}\thanksref{T2,m1}}
\and
\author[B]{\fnms{Daniel}~\snm{Remenik}\corref{}\ead[label=e3]{dremenik@dim.uchile.cl}\thanksref{T3,m2}}
\runauthor{J.~Ortmann, J.~Quastel and D.~Remenik}
\affiliation{University of Toronto\thanksmark{m1} and Universidad de
Chile\thanksmark{m2}}
\thankstext{T1}{Supported in part by the
Natural Sciences and Engineering Research Council of Canada.}
\thankstext{T2}{Supported by the Natural Sciences and Engineering Research
Council of Canada,
the I. W. Killam Foundation and the Institute for Advanced Study.}
\thankstext{T3}{Supported in part by Fondecyt Grant 1120309,
by Conicyt Basal-CMM and by
Programa Iniciativa
Cient\'ifica Milenio grant number NC130062 through Nucleus Millenium
Stochastic Models of
Complex and Disordered Systems.}
\address[A]{J.~Ortmann\\
J.~Quastel\\
Department of Mathematics\\
University of Toronto\\
40 St. George Street\\
Toronto, Ontario M5S 2E4\\
Canada\\
\printead{e1}\\
\phantom{E-mail:\ }\printead*{e2}}
\address[B]{D.~Remenik\\
Departamento de Ingenier\'ia Matem\'atica\\
\quad and Center for Mathematical Modeling\\
Universidad de Chile\\
Av. Blanco Encala\-da 2120\\
Santiago\\
Chile\\
\printead{e3}}
\end{aug}

%
\received{\smonth{10} \syear{2014}}
%
\revised{\smonth{1} \syear{2015}}

%
\begin{abstract}
We obtain exact formulas for moments and generating functions of the
height function of
the asymmetric simple exclusion process at one spatial point, starting
from special
initial data in which every positive even site is initially occupied.
These complement
earlier formulas of E. Lee [\textit{J. Stat. Phys.} \textbf{140} (2010) 635--647]
but, unlike those formulas, ours
are suitable in
principle for asymptotics. We also explain how our formulas are
related to divergent
series formulas for half-flat KPZ of Le Doussal and Calabrese
[\textit{J. Stat. Mech.} \textbf{2012} (2012) P06001], which we
also recover using the methods of this paper. These generating
functions are given as a series without any apparent Fredholm
determinant or Pfaffian structure. In the long time limit, formal asymptotics
show that the fluctuations are given by the Airy$_{2\to1}$ marginals.
\end{abstract}

%
\begin{keyword}[class=AMS]
\kwd{60K35}
\kwd{82C22}
\kwd{82C23}
\kwd{82B23}
\kwd{60H15}
\end{keyword}
\begin{keyword}
\kwd{Interacting particle systems}
\kwd{Kardar--Parisi--Zhang universality class}
\kwd{asymmetric simple exclusion process}
\kwd{flat initial data}
\end{keyword}
\end{frontmatter}

\section{Introduction}

\label{sec:intro}

The one-dimensional asymmetric simple exclusion process (ASEP) is a
continuous time Markov
process with state space $\{0,1\}^{\zz}$, the $1$'s being thought of
as particles and the $0$'s as holes. Each particle has an independent
exponential clock
which rings at rate one. When it rings, the particle chooses to attempt
to jump one site
to the right with probability $p\in[0,1]$, or one site to the left with
probability $q=1-p$.
However, the jump is only executed if the target site is empty;
otherwise, the jump is
suppressed and the particle must wait for the alarm to ring again. If
$q=1$, $p=0$ (or
$q=0$, $p=1$, but we will assume for convenience that $q\ge p$), it is
called the totally
asymmetric simple exclusion process (TASEP); if $0<q\neq p$ it is the
(partially)
asymmetric simple exclusion process (ASEP); if $q=p=1/2$ it is the
symmetric simple
exclusion process (SSEP). We denote by $\eta_t(x)=1$ or $0$ the
presence or absence of a
particle at $x\in\zz$ at time $t$. The state of the system is
completely determined at time $t>0$ by the
initial data $\eta_x(0)$, $x\in\zz$, together with the family of
exponential clocks; for
more details on the construction of the process, we refer the reader to
\cite{ligg1}.
Given $\eta\in\{0,1\}^{\zz}$, we define $\hat\eta\in\{-1,1\}^{\zz
}$ by
$\hat{\eta}(x)=2\eta(x)-1$. The \emph{height function} of ASEP is
defined in terms of $\hat\eta_t$ by
%
%
\begin{equation}
\label{eq:defofheight} h(t,x) = %
\cases{\displaystyle 2N^{\mathrm{ flux}}_0(t) +
\sum_{0<y\leq x}\hat{\eta}_t(y), &\quad  $x>0$,
\vspace*{2pt}
\cr
2N^\mathrm{ flux}_0(t), &\quad $ x=0$,\vspace*{2pt}
\cr
\displaystyle 2N^\mathrm{ flux}_0(t)- \sum_{x<y\leq0}
\hat{\eta}_t(y), &\quad  $x<0$,} %
\end{equation}
where $N^\mathrm{ flux}_0(t)$ is the net number of particles which crossed
from site $1$ to $ 0$ up to
time $t$, meaning that particle jumps $1\to0$ are counted as $+1$ and
jumps $0\to1$ are
counted as $-1$.

ASEP is an important member of the one-dimensional
Kardar--Parisi--Zhang (KPZ) universality
class. This is a broad class of one-dimensional driven diffusive
systems, or
stochastic growth models, characterized by unusual, but universal
asymptotic fluctuations.
These should be of size $t^{1/3}$, and decorrelate on spatial scales of
$t^{2/3}$, with
special distributions in the long time limit, usually given in terms of
Fredholm determinants, which only
depend on the initial data class. There are a few special classes of
initial data
characterized by scale invariance: \emph{curved} (or \emph{step}),
corresponding to
starting with particles at every nonnegative site; \emph{flat} (or
\emph{periodic}),
corresponding to starting with particles at all even sites; and \emph
{stationary},
corresponding to starting with a product Bernoulli measure. In
addition, there are three
crossover classes: \emph{curved${}\to{}$flat}, \emph{curved$\to
$stationary} and
\emph{flat${}\to{}$stationary}; corresponding to putting two different
initial conditions on
either side of the origin. Based on exact computations for TASEP and a
few other models with
special determinantal (Schur) structure, the asymptotic spatial
fluctuations in all six
cases are known to be given by the Airy processes, a family of
processes defined through
their finite dimensional distributions which are given by specific Fredholm
determinants. The full space--time limit in this KPZ scaling $\ep^{1/2}
h(\ep^{-3/2} t,
\ep^{-1} x)$ is believed to be a Markov process known as the \textit{KPZ
fixed point}. For
more details, see the reviews \cite{quastelCDM,corwinReview,quastelRem-review}.

Within the universality class, the KPZ equation
\[
\partial_t h= \frac{1}2\partial_x^2
h +\frac{\gamma}2(\partial_xh)^2 + \xi,
\]
where $\xi$ is space--time white noise, plays a special role as a
heteroclinic orbit
connecting the Edwards--Wilkinson (linear) fixed point $\partial_t h=
\frac{1}2\partial_x^2
h + \xi$ to the (nonlinear and poorly understood) KPZ fixed point. It
can be obtained
from other models with adjustable nonlinearity or noise in the diffusive
($t=\ep^{-2}T$, $x=\ep^{-1}X$) weakly asymmetric, or weak noise
limit, with
rigorous proofs available in a few cases
\cite{bertiniCancrini,berGiaco,acq,akq2,hairer,demboLi,hairerReg,mqrScaling}.

The importance of ASEP in this context is that it has an adjustable
nonlinearity
\[
\gamma= q-p.
\]
Although in the case $\gamma>0$ it does not have a determinantal
structure, somewhat
surprisingly exact formulas have been discovered for the distribution
of the height
function of ASEP at a fixed $t>0$ and $x\in\zz$ for certain initial
data, starting with
the work of Tracy and Widom in 2008 \cite
{tracyWidomASEP2,tracyWidomASEP1}. The first
formula was for the \emph{step} case $\eta_0^\mathrm{ step} (x) = \mathbf
{1}_{x\in \zz_{>0}}$. In the
weakly asymmetric limit exact formulas were obtained for the one-point
distribution of the
KPZ equation with so called \emph{narrow wedge} initial data
(corresponding to the curved
class); see \cite{acq} and also \cite{sasamSpohn}. In the $t\to
\infty$
limit, one obtains
the Tracy--Widom GUE distribution. An analogous procedure was then
performed on the step
Bernoulli, or curved${}\to{}$stationary case for ASEP, corresponding to
half-Brownian initial
data for KPZ; the $t\to\infty$ limit in this case gives the
Airy$_{2\to
{\mathrm{ BM}}}$
marginals, or BBP transitional distributions~\cite
{corwinQuastelRarefaction}. Parallel
computations were performed on the physics side using the nonrigorous
replica method. The
case of Brownian initial data for KPZ (corresponding to stationary
ASEP) has also recently
been completed in the physics \cite{imamSasamStationary} and mathematics
\cite{borCorFerrVeto} literatures. It should be
emphasized that these are formulas for one point distributions only,
and for very special
initial data. So far, multipoint distributions have resisted rigorous
analysis, though
some nonrigorous attempts have been made
\cite{prolhacSpohn,dotsenko2point,dotsenko2time}.

Among the primary scaling invariant initial data at the KPZ level, this
left the flat and
half-flat cases. In \cite{cal-led,leDoussalHF}, Le Doussal and
Calabrese gave formulas
for the one point height distribution of KPZ for the half-flat and flat
initial data via
the replica method. Their half-flat formula is an uncontrolled
divergent series, with no apparent Fredholm
structure. As such, it is a pure formalism, and is mainly used as an
intermediate step in order to
obtain a Fredholm Pfaffian formula for the flat initial condition, by
scaling the wedge to infinity, that is, looking farther and farther
into the flat
region.

Here, we will work directly
with ASEP, which in particular can be regarded as a microscopic version
of KPZ
\cite{berGiaco}, and where one can avoid the problems associated with
the nonsummable
moments. Later, in Section~\ref{sec:app-bose}, we will discuss how the
methods we will use
can be applied in the case of KPZ, yielding some of the formulas
appearing in
\cite{cal-led,leDoussalHF}.

We will be primarily concerned with the \emph{half-flat} initial condition,
%
%
\begin{equation}
\label{eq:halfflat-initial-condition} \eta^{\mathrm{ h\mbox{-}fl}}_0(x)=\mathbf{1}_{x\in2\zz_{>0}}.
\end{equation}
The superscript h\mbox{-}fl will be used for probabilities and expectations
computed with respect to
this initial condition. The limit to the flat initial condition $\eta
_0^\mathrm{
flat}(x)=\mathbf{1}_{x\in2\zz}$ will be pursued in an upcoming paper.

E. Lee's thesis already contains exact formulas
for the quantities we are interested in. Here, and in the rest of the
paper, we set
\[
\tau=\frac{p}{q}\in(0,1).
\]
%

\begin{thmm}[(\cite{lee})]
%
%
\begin{eqnarray}
\label{eq:final2}&& \mathbb{P}^{\mathrm{ h\mbox{-}fl}}\bigl(h(t,0) \ge2m-x\bigr)
\nonumber\\
&&\qquad=
(-1)^m\sum_{k \geq
m}\frac{\tau^{(k-m)(k-m+1)/2}}{(1+\tau)^{k(k-1)}k!}
\pmatrix{k-1
\cr
k-m} _\tau
\nonumber
\\[-8pt]
\\[-8pt]
\nonumber
&&\qquad\quad{}\times\int_{C_R^k} \prod_{i\neq j}
\frac{\xi_j-\xi_i}{p+q\xi_i\xi_j
-\xi_i}\prod_{i}\frac{\xi_i^{x}e^{t\varepsilon
(\xi_i)}}{(1-\xi_i)(\xi_i^2-\tau)}\\
&&\hspace*{28pt}\qquad\quad{}\times\prod
_{i<j}\frac{1+\tau- (\xi
_i +\xi
_j)}{\tau
-\xi_i\xi_j}\prod
_i d\xi_i,
\nonumber
\end{eqnarray}
where
%
%
\begin{equation}
\varepsilon(\xi_i):= p\xi_i^{-1} + q
\xi_i -1.\label{eq:vareps}
\end{equation}
$C_{R}$ is a contour large enough to contain all the poles of the
integrand, and ${n\choose k} _\tau=\frac{n_\tau!}{k_\tau
!(n-k)_\tau!}$ with
the $\tau$-factorial
$n_\tau!$ defined in \eqref{eq:qfact}.
\end{thmm}

These formulas are similar in structure to earlier formulas of \cite
{tracyWidomASEP1}.
However, such formulas turn out not to be conducive to asymptotics
analysis. They
need considerable ``postproduction'' before the asymptotic behaviour
can be
extracted \cite{tracyWidomASEP2,tracyWidomASEP3}, and no one
has been able to figure out how to do this for \eqref{eq:final2}, nor to
extract the relevant asymptotics (even formally).

Our main result is an explicit formula for the one-point distribution
in the half-flat
case, expressed as a certain series which has a structure reminiscent
of a Fredholm
determinant. In an upcoming paper, we will use these formulas to obtain
analogous moment
formulas in the flat case and, furthermore, a Fredholm Pfaffian
formula for a certain
transform of the height function. Formal asymptotics lead to the expected
results in the $t\to\infty$ and weakly asymmetric limits, but turning
them into rigorous proofs presents some
considerable technical challenges and is left for future work (see
the \hyperref[sec:app-2to1]{Appendix} for a discussion of the large time case).

Formulas for the half-flat case can be obtained by the method of \cite
{bcs}, together
with an ansatz coming from a study of the mechanics of \eqref
{eq:final2}. Let
%
%
\begin{equation}
\label{eq:Nx} N_x(t)=\sum_{y=-\infty}^x
\eta_t(y)
\end{equation}
be the total number of particles to the left of $x$ at time $t$. It is
not hard to check
that when all particles start to the right of the origin, $N^\mathrm{
flux}_0(t)=N_0(t)$, and
thus by~\eqref{eq:defofheight}
%
%
\begin{equation}
\label{eq:NflN0} h(t,x)=2N_x(t)-x
\end{equation}
in the half-flat case. Define
\[
\label{eq:defQs} \wt Q_x(t)=\frac{\tau^{N_x(t)}-\tau^{N_{x-1}(t)}}{\tau-1}=\tau ^{N_{x-1}(t)}
\eta_x(t).
\]

%
\begin{thmm}
\label{Thm:FormulaU}
Consider ASEP with half-flat initial condition as in
\eqref{eq:halfflat-initial-condition}. Then for any $\vec x\in\zz^k$
we have
%
%
\begin{eqnarray}
\label{Eq:MomentsQTilde} &&\ee^{\mathrm{ h\mbox{-}fl}} \bigl[\wt{Q}_{x_1}(t)\cdots
\wt{Q}_{x_k}(t) \bigr]
\nonumber
\\
&&\qquad=\frac
{\tau^{({1}/{2})k(k-1)}}{(2\pi\I)^{k}}\\
&&\qquad\quad{}\times\int_{C_{1,\rho}^k}
\,d\vec z \prod_{1\leq a<b\leq k} \frac{z_a- z_b}{z_a - \tau z_b }
\frac{1-z_az_b}{1-\tau z_az_b} \prod_{a=1}^k
\frac{1}{\tau z_a^2-1}f_{x_a,t} (z_a),\nonumber
\end{eqnarray}
where $C_{1,\rho}$ is a circle around 1 with radius
$0<\rho<\min\{\tau^{-1/2}-1,(1+\tau)^{-1}\}$, $C_{1,\rho}^k$ denotes
the product of $k$
copies of $C_{1,\rho}$,
\[
\label{Eq:DefFxt} f_{x,t}(z) = \biggl( \frac{1-\tau z}{ 1-z}
\biggr)^{x-1} e^{\wt\ep
(z) t},
\]
and $\wt\varepsilon$ is defined in terms of the function $\varepsilon$
given in
\eqref{eq:vareps} by
%
%
\begin{equation}
\label{Eq:DefEpTilde} \wt\varepsilon(z) = \varepsilon \biggl(\frac{1-\tau z}{1-z} \biggr)
= p\frac{1-z}{1-\tau z}+q\frac{1-\tau z}{1-z}-1.
\end{equation}
\end{thmm}

For simplicity, throughout the rest of the paper we will omit the bound
on the indices in
products such as $\prod_{1\leq a\leq k}$ and $\prod_{1\leq a<b\leq k}$
when no confusion can arise
and the factors involved in the products are defined in terms of a
collection of $k$
variables. A similar convention will sometimes be used for sums.
Additionally, we will
continue using the notation $C^k$ for the product of $k$ copies of a given
contour $C$ in the complex plane.

An analogous formula holds for the stochastic heat equation/KPZ/delta
Bose gas; see
Section~\ref{sec:app-bose} for details. On the other hand, the
analogous ansatz does not
work in the case of $q$-TASEP and the O'Connell--Yor semi-discrete
polymer, at least with the most straightforward candidates for half-flat
initial conditions.

The interesting new term here over earlier formulas \cite{bcs,borCor}
is $\prod_{a<b
}\frac{1-z_az_b}{1-\tau z_az_b}$, which together with the factor
$\prod_a \frac{1}{\tau z_a^2-1}$ allows us to recover the periodic initial
data. This term
leads to the
double product $\prod_{a<b}\mfh(w_a,w_b;s_a,s_b)$ appearing in \eqref
{eq:nuk-intro}
below, which is the obstacle to making long time limit fluctuations
rigorous (see
the \hyperref[sec:app-2to1]{Appendix}). Factors of this form were fortuitously
absent from earlier
formulas for step and step Bernoulli initial data, which only
contained the double
product $\prod_{a<b}\frac{z_a-z_b}{z_a-\tau z_b}$; this last factor
turns into the
determinant in \eqref{eq:nuk-intro} and this makes it much easier to deal
with. Similar expressions have also proved to be an obstacle in the
replica formulas
\cite{cal-led,leDoussalHF} for half-flat and flat initial data, as
well as for expressions for
multipoint distributions \cite{dotsenko2time,dotsenko2point,dotsenkoNpoint}.\vspace*{1pt}

The formula for $\ee^{\mathrm{ h\mbox{-}fl}} [\wt{Q}_{x_1}(t)\cdots
\wt{Q}_{x_k}(t)
 ]$ can be
used to write a formula for the moments of $\tau^{N_x(t)}$ by using
ideas of
\cite{imamSasamASEPduality,bcs}. The result is given in Section~\ref{sec:mom} as
Proposition~\ref{prop:NestedPlusSmall}. The formula for $\ee[\tau
^{kN_x(t)}]$ is given as
a nested contour integral (see Figure~\ref{fig:nested}). As given, such
a formula is
suitable neither for asymptotic analysis (not even at a formal level)
nor for our later
goal of deriving a formula for the full flat case. In order to obtain a
formula where all
the contours coincide we will expand the nested contours so that they
all coincide with
largest one. The resulting formula amounts to computing the residue
expansion associated
to the poles that we cross as we perform this deformation. It is given
in Proposition~\ref{prop:muk} as a sum of multiple contour integrals indexed by
partitions. After some
rewriting, this formula leads to our main result for ASEP with
half-flat initial
data. Define the following functions:
%
%
\begin{eqnarray}
\label{eq:germans-hf} %
\mff(w;n)&=&(1-\tau)^{n}e^{(q-p)t [{1}/{(1+w)}-{1}/{(1+\tau
^{n}w}) ]}
\biggl(\frac{1+\tau^{n}w}{1+w} \biggr)^{x-1},\nonumber
\\
\mfg(w;n)&=&\frac{ (-w;\tau )_\infty}{ (-\tau
^{n}w;\tau )_\infty} \frac{ (\tau^{2n}w^2;\tau )_\infty}{ (\tau
^{n}w^2;\tau )_\infty},
\\
\mfh(w_1,w_2;n_1,n_2)&=&
\frac{ (w_1w_2;\tau )_\infty
(\tau ^{n_1+n_2}w_1w_2;\tau )_\infty}{
 (\tau^{n_1}w_1w_2;\tau )_\infty (\tau
^{n_2}w_1w_2;\tau )_\infty},\nonumber
\end{eqnarray}
where the infinite \emph{$q$-Pochhammer symbols} are defined as
\[
\label{Eq:DefQPoch} (a;q )_\infty=\prod_{n=0}^\infty
\bigl(1-q^na\bigr).
\]
Note that $\mfg$ and $\mfh$ can be written in
terms of ratios of finite $q$-Pochhammer symbols, but it will more
convenient for us to
write them in this form. The formulas for $\mfg$ and $\mfh$ can
alternatively be written as ratios of \emph{$q$-Gamma functions},
\[
\label{Eq:DefQGamma} \Gamma_q(x) = \frac{(1-q)^{1-x}  (q;q )_\infty}{
(q^x;q )_\infty},
\]
which converge (uniformly on compact sets) to the usual Gamma function
as \mbox{$q\to1$}. We also define the \emph{$q$-factorial}
%
%
\begin{equation}
\label{eq:qfact} m_q!=\frac{\prod_{a=1}^k(1-q^a)}{(1-q)^k}.
\end{equation}
For later use, we further introduce the \emph{$q$-exponential function}
%
%
\begin{equation}
\label{eq:qexp} e_q(x)=\frac{1}{ ((1-q)x;q )_\infty}=\sum
_{k=0}^\infty \frac{x^k}{k_q!},
\end{equation}
where the second equality only holds for $|x|<1$ and amounts to the
$q$-Binomial theorem
(see, e.g., Theorem~10.2.1 in \cite{AndrewsAskeyRoy}). As $q\to
1$, this function
converges to the usual exponential function, uniformly on $(-\infty,A]$
for any $A$. In
keeping with the standard usage we have used the parameter $q$ in the
definition of these
$q$-deformed functions, but in all that follows the parameter $\tau$
will appear in place
of $q$.

%
\begin{thmm}\label{thmm:main-hf}
Consider ASEP with half-flat initial condition as in
\eqref{eq:halfflat-initial-condition} and let $m\in\zz_{\geq0}$. Then
%
%
\begin{equation}
\label{eq:stpt1} \ee^{\mathrm{ h\mbox{-}fl}} \bigl[\tau^{mN_{x}(t)} \bigr]
=m_\tau!\sum_{k=0}^{m}
\nnu^{\mathrm{ h\mbox{-}fl}}_{k,m}(t,x)
\end{equation}
with
\begin{eqnarray}
\label{eq:nuk-intro} \nnu^{\mathrm{ h\mbox{-}fl}}_{k,m}(t,x)&=&\frac{1}{k!}
\mathop{\sum_{n_1,\ldots,n_k\geq1}}_{
n_1+\cdots+n_k=m}\frac{1}{(2\pi\I)^{k}}
\int_{\gamma_{-1,0}^k}\,d\vec w\, \det \biggl[\frac{-1}{w_a\tau
^{n_a}-w_b}
\biggr]_{a,b=1}^{k}
\nonumber
\\[-8pt]
\\[-8pt]
\nonumber
&&{}\times\prod_a\mff(w_a;n_a)
\mfg(w_a;n_a)\prod_{a<b}
\mfh(w_a,w_b;n_a,n_b),
\nonumber
\end{eqnarray}
where $\gamma_{-1,0}$ is a (positively oriented) contour around $-1$
and $0$, strictly
contained inside the circle of radius $\tau^{-1/2}$, which does not
include any other
singularities of the integrand.
\end{thmm}

The contour $\gamma_{-1,0}$ in the theorem can for example be chosen to
be a circle around
the origin with radius in $(1,\tau^{-1/2})$. In fact, the determinant
clearly never
vanishes for this choice, and one can check that all the other
singularities of the
integrand, except for $w_a=0$ and $w_a=-1$, are outside this contour.

With a formula for the moments of $\tau^{N_{x}(t)}$ at our disposal we
are ready to form a
generating function, namely the \emph{$\tau$-Laplace transform} of
$\tau
^{N_x(t)}$. The
formula involves a Mellin--Barnes integral representation of the
infinite sums in
$n_1,\ldots,n_k$ appearing in \eqref{eq:nuk-intro} after summing over
$m\geq0$.

%
\begin{thmm}\label{thmm:gen-hf}
Let $\zeta\in\cc\setminus\rr_{\geq0}$. Then, for $e_\tau$ as in
\eqref{eq:qexp},
%
%
\begin{eqnarray}
\label{eq:gen-hf} &&\ee^{\mathrm{ h\mbox{-}fl}} \bigl[e_\tau\bigl(\zeta
\tau^{N_x(t)}\bigr) \bigr]\nonumber\\
&&\qquad =\sum_{k=0}^{\infty}
\frac{1}{k!}\frac{1}{(2\pi\I)^{2k}}\int_{({1}/{2}+\I\rr
)^k}\,d\vec s \int
_{\gamma_{-1,0}^k}\,d\vec w\, \det \biggl[\frac{-1}{w_a\tau
^{s_a}-w_b}
\biggr]_{a,b=1}^{k}
\\
&&\qquad\quad{}\times\prod_a(-\zeta)^{s_a}
\mff(w_a;s_a)\mfg(w_a;s_a)\prod
_{a<b}\mfh (w_a,w_b;s_a,s_b).
\nonumber
\end{eqnarray}
\end{thmm}

Now set $\zeta=-\tau^{-t/4-t^{2/3}x/2+t^{1/3}r/2}$. Since $e_\tau
(z)\longrightarrow0$ as
$z\to-\infty$ and $e_\tau(z)\longrightarrow1$ as $z\to0$ for fixed
$\tau
$, uniformly in
$z\in(-\infty,0]$ we have (see \cite{borCor}, Lemma~4.1.39)
\begin{eqnarray}\label{eq:eepp-hf}\quad
&&\lim_{t\to\infty}\ee^{\mathrm{ h\mbox{-}fl}} \bigl[e_\tau\bigl(-
\tau ^{N_{t^{2/3}x}(t/\gamma
)-({1}/{4})t-({1}/{2})t^{2/3}x+({1}/{2})t^{1/3}r-
({1}/{4})t^{1/3}x^2\mathbf{1}_{x\leq 0}}\bigr) \bigr]
\nonumber
\\[-8pt]
\\[-8pt]
\nonumber
&&\qquad=\lim_{t\to\infty}\pp^{\mathrm{ h\mbox{-}fl}}
\biggl(\frac{h(t/\gamma
,t^{2/3}x)-({1}/2)t-({1}/{2})t^{1/3}x^2\mathbf{1}_{x\leq0}}{t^{1/3}}\geq-r
\biggr),
\nonumber
\end{eqnarray}
where we recall $\gamma=q-p$. In the \hyperref[sec:app-2to1]{Appendix}, we show
that a formal steepest
descent analysis of the right-hand side of \eqref{eq:gen-hf} gives (a
scaled version of)
the one-point marginals of the Airy$_{2\to1}$ process ${\mathcal
{A}}_{2\to1}(x)$.

\textit{Outline}. The rest of the paper is organized as follows.
Section~\ref{sec:ansatz} contains the proof of Theorem~\ref
{Thm:FormulaU}. In Section~\ref{sec:mom} we will use the formula
obtained in Theorem~\ref
{Thm:FormulaU} to derive the
moment formula given in Theorem~\ref{thmm:main-hf}, while in
Section~\ref{sec:genfn} we
will derive the formula for the $\tau$-Laplace transform of $\tau
^{N_x(t)}$ (Theorem~\ref{thmm:gen-hf}). Section~\ref{sec:app-bose} explains how the methods
used for ASEP can
be applied to the case of the SHE/KPZ equation (or, more precisely, the
delta Bose gas)
and discusses the relation with the work of Le Doussal and Calabrese.
Finally, the \hyperref[sec:app-2to1]{Appendix} contains the formal derivation of
the limiting
fluctuations for ASEP
with half-flat initial condition.

\section{Contour integral ansatz}\label{sec:ansatz}

To prove Theorem~\ref{Thm:FormulaU}, we will use Proposition~4.10 of
\cite{bcs}, which
shows that $\ee^{\mathrm{ h\mbox{-}fl}} [\wt{Q}_{x_1}(t)\cdots\wt{Q}_{x_k}(t)
]$ can be
represented as the solution of a certain evolution equation with
boundary conditions. We
describe this result next.

Let $\eta_0$ be an ASEP configuration with a leftmost particle and
consider ASEP started
with $\eta_0$ as initial condition. Let $\tilde u_0(\vec x)=\prod_{a=1}^k\tau^{N_{x_a-1}(0)}\eta_{x_a}(0)$ (where, of course, $N_x(0)$
is computed
with respect to the initial condition $\eta_0$). Consider the following
system of
differential equations:
\begin{longlist}[(1)]
\item[(1)]
For all $\vec x\in\zz^k$ and $t\geq0$, writing
$\vec{x}^\pm_\ell=(x_1,\ldots,x_\ell\pm1,\ldots,x_k)$,
\[
\frac{d}{d t} \tilde u(t,\vec x) = \sum_{j=1}^k
\bigl[p \tilde u\bigl(t,\vec x_j^{-}\bigr) + q \tilde u
\bigl(t,\vec{x}^+_j\bigr)-\tilde u(t,\vec x) \bigr].
\]
\item[(2)]
For all $\vec x\in\zz^k$ such that there
exists $\ell<k$ with
$x_{\ell+1}=x_\ell+1$,
\[
p\tilde u\bigl(t,\vec{x}^-_{\ell+1}\bigr) + q \tilde u\bigl(t,
\vec{x}^+_\ell\bigr) = \tilde u(t,\vec{x}).
\]
\item[(3)]
There exist constants $c,C,\delta>0$ such
that for all $\vec
x\in\zz^k$ with $x_1<x_2<\cdots<x_k$ and $t\in[0,\delta]$,
\[
\bigl\llvert \tilde u(t,\vec x)\bigr\rrvert \leq C e^{c\sum_{j} \llvert x_a\rrvert }.
\]
\item[(4)]
For all $\vec x\in\zz^k$ such that
$x_1<x_2<\cdots<x_k$ we have
\[
\tilde u(0,\vec{x})=\tilde u_0(\vec x).
\]
\end{longlist}

%
\begin{prop}[(\cite{bcs})]\label{prop:bcs}
Suppose that $\tilde u(t,\vec x)$ solves \emph{(1)--(4)}. Then for all
$\vec{x}\in\zz^k$ such
that $x_1<x_2<\cdots<x_k$ we have
\[
\ee^{\eta_0} \bigl[\wt{Q}_{x_1}(t)\cdots\wt{Q}_{x_k}(t)
\bigr]=\tilde u(t,\vec x),
\]
where the superscript on the left-hand side means that ASEP is started
with initial
condition $\eta_0$.
\end{prop}

We proceed now to the proof of our formula for
$\ee^{\mathrm{ h\mbox{-}fl}} [\wt{Q}_{x_1}(t)\cdots\wt{Q}_{x_k}(t)
 ]$.

\begin{pf*}{Proof of Theorem~\ref{Thm:FormulaU}}
In view of Proposition~\ref{prop:bcs} and \eqref{Eq:MomentsQTilde}, we
need to check that
%
%
\begin{eqnarray}\label{eq:utilde}
\tilde u(t;\vec x)&:=&\frac{\tau^{({1}/{2})k(k-1)}}{(2\pi\I
)^{k}}\int_{C_{1,\rho}^k} \,d\vec z
\prod_{1\leq a<b\leq k} \frac{z_a- z_b}{z_a - \tau z_b } \frac{1-z_az_b}{1-\tau z_az_b}
\nonumber
\\[-8pt]
\\[-8pt]
\nonumber
&&{}\times\prod_{a=1}^k \frac{1}{\tau z_a^2-1}f_{x_a,t}
(z_a)
\end{eqnarray}
satisfies (1)--(4) with $\tilde u_0$
defined in terms of
the half-flat initial condition $\eta_0(x)=\mathbf{1}_{x\in2\zz
_{>0}}$. A
straightforward
computation shows that in this case
%
%
\begin{equation}
\label{eq:tildeu0-hf} \tilde u_0(\vec x)=\prod
_{a=1}^k \mathbf{1}_{x_a\in2\zz_{>0} }
\tau^{\sum_{y=-\infty}^{x_a-1} \eta_y(0)} = \tau^{-k}\prod_{a=1}^k
\mathbf{1}_{x_a\in2\zz_{>0} } \tau^{ ({1}/2) x_a}.
\end{equation}
We will denote the integrand in \eqref{eq:utilde} by $I_{k,t}(\vec
{x};\vec{z})$,
that is,
%
%
\begin{equation}
I_{k,t}(\vec{x};\vec{z})=\prod_{a<b}
\frac{z_a- z_b}{z_a - \tau z_b } \frac{1-z_az_b}{1-\tau z_az_b} \prod_a
\frac{1}{\tau z_a^2-1} f_{x_a,t} (z_a).\label{eq:Ikt}
\end{equation}
Additionally, we will write $\vec{x}^{(i_1,\ldots,i_\ell)}$ and
$\vec
z^{(i_1,\ldots,i_\ell)}$
to denote, respectively, the vectors $\vec x$ and $\vec{z}$ with the
components $i_1,\ldots,i_\ell$ removed.

Computing $\frac{d}{dt} \tilde u(t,\vec{x})$ introduces a factor
$\sum_{\ell=1}^k\wt\ep(z_\ell)$ in front
of the integrand. Similarly, computing $\tilde u(t,\vec{x}^\pm_\ell)$
introduces a factor $ (\frac{1-\tau z_\ell}{1-z_\ell}
)^{\pm
1}$ in front of the
integrand. Hence, (1) is satisfied if we can show that
\[
\sum_{\ell=1}^k\wt\ep(z_\ell) =
\sum_{\ell=1}^k \biggl[p\frac{1-z_\ell}{1-\tau z_\ell}+q
\frac
{1-\tau
z_\ell}{1-z_\ell}-1 \biggr].
\]
But this follows immediately from the definition of
$\wt\ep$; see \eqref{Eq:DefEpTilde}

For (2), let $\vec{x}\in\zz^k$ and suppose that there
exists $\ell$ such
that $x_{\ell+1}=x_\ell+1$. Then using the above computation of
$\tilde
u(t,\vec{x}^\pm_\ell)$, we have
%
%
\begin{eqnarray}
\label{Eq:Antisymm}&& p \tilde u\bigl(t,\vec{x}^{-}_{\ell+1}\bigr) +q
\tilde u\bigl(t,\vec{x}^+_{\ell}\bigr) -\tilde u(t,\vec{x})\nonumber\\
&&\qquad =
\frac{\tau^{({1}/{2})k(k-1)}}{(2\pi\I)^{k}}\int_{C_{1,\rho}^k} \,d\vec{z}\, I_{k,t}\bigl(
\vec{x}^-_{\ell+1};\vec z\bigr)
\\
&&\qquad\quad{}\times \biggl[p+q\frac{1-\tau z_\ell}{1-z_\ell}\frac{1-\tau z_{\ell
+1}}{1-z_{\ell+1}}- \frac{1-\tau z_{\ell+1}}{1-z_{\ell+1}}
\biggr].\nonumber
\end{eqnarray}
We need to show that the integral vanishes. The expression inside the
brackets equals
$\frac{(q-p)(z_\ell-\tau z_{\ell+1})}{(1-z_\ell)(1-z_{\ell+1})}$. Note
that the factor
$z_\ell-\tau z_{\ell+1}$ cancels a like factor in the denominator of
the product
$\prod_{a<b}\frac{z_a-z_b}{z_a-\tau z_b}$\vspace*{1pt} coming from $I_{k,t}(\vec
{x}^-_{\ell+1};\vec
z)$, and thus (using the fact that $x_{\ell+1}=x_\ell+1$) the
integrand in
\eqref{Eq:Antisymm} can be rewritten as
\[
\frac{(q-p)(z_\ell-z_{\ell+1})(1-z_\ell z_{\ell+1})}{(1-z_\ell
)(1-z_{\ell+1})(1-\tau
z_\ell z_{\ell+1})}f_{x_\ell,t}(z_\ell)f_{x_\ell,t}(z_{\ell
+1})G
\bigl(\vec {x}^{(\ell,\ell+1)},\vec{z}^{(\ell,\ell+1)}\bigr),
\]
where, as suggested by the notation, the factor
$G(\vec{x}^{(\ell,\ell+1)},\vec{z}^{(\ell,\ell+1)})$ does not depend
on $x_\ell$,
$x_{\ell+1}$, $z_\ell$ and $z_{\ell+1}$. This expression is
antisymmetric in $z_\ell,z_{\ell+1}$, and thus its integral over
$(z_\ell,z_{\ell+1})\in C_{1,\rho}^2$ must vanish. This shows that the
integral in~\eqref{Eq:Antisymm} is zero, proving (2).

(3) follows directly from the form of $f_{x,t}$ and the
facts that
$C_{1,\rho}$ is compact and that the integrand
is continuous in $\vec z\in C_{1,\rho}^k$.

We turn now to (4). Note that when $t=0$ the essential
singularity in
the exponent of $f_{x,t}$ in $I_{k,t}$ disappears [see \eqref
{eq:Ikt}], and thus we can
evaluate the integral by computing residues.

First, if $x_1\leq1$ then $f_{x_1,0}(z_1)$ has no pole at $z_1=1$.
Hence, the integrand
is analytic in $z_1$ inside $C_{1,\rho}$, and thus the integral is 0. Since
$x_1<\cdots<x_k$, this accounts for the condition that all $x_a$'s be
at least 2. So let
us assume now that $2\leq x_1<\cdots<x_k$. We will evaluate the $z_k$
integral first, by
expanding the contour to infinity. Note that, thanks to the decay
coming from the factor
$(\tau z_k^2-1)^{-1}$ there is no pole at infinity, and thus the
integral equals minus
the sum of the residues of the poles of the integrand outside
$C_{1,\rho}$.

In $z_k$, the poles are $\pm\tau^{-1/2}$, $\tau^{-1} z_\ell$ and
$\tau
^{-1} z_\ell^{-1}$
for $\ell<k$. The condition imposed on $\rho$ implies that all these
poles lie outside
the contour. Consider first the poles at $z_k=\tau^{-1}z_\ell$, $\ell<k$.
The residue of $I_{k,0}$ at
this point is given by
\begin{eqnarray*}
&&I_{k-1,0}\bigl(\vec{x}^{(k)};\vec{z}^{ (k)}\bigr)
\mathop{\prod_{a<k}}_{ a \ne
\ell} \biggl(
\frac{z_a-\tau^{-1} z_\ell}{z_a -
z_\ell}\frac{1-\tau^{-1}z_az_\ell}{1-z_az_\ell} \biggr)
\frac{(1-\tau)z_\ell(1-\tau^{-1}z_\ell^2)}{1-z_\ell^2}\\
&&\hspace*{78pt}\quad{}\times \frac
{
({(1-
z_\ell)}/{(1-\tau^{-1}z_\ell)} )^{x_k-1}}{\tau^{-1}z_\ell^2-1}
\\
&&\qquad=I_{k-2,0}\bigl(\vec{x}^{( \ell, k)};\vec{z}^{(\ell,k)}\bigr)\prod
_{a=1}^{\ell
-1} \biggl(\frac{z_a-\tau^{-1} z_\ell}{z_a -
\tau z_\ell}
\frac{1-\tau^{-1}z_az_\ell}{1-\tau z_az_\ell} \biggr) \\
&&\qquad\quad{}\times \prod_{b=\ell+1}^{k-1}
\biggl(\frac{\tau^{-1} z_\ell-z_b}{z_\ell-
\tau z_b}\frac{1-\tau^{-1}z_bz_\ell}{1-\tau z_bz_\ell} \biggr)
\\
&& \qquad\quad{}\times\frac{(1-\tau)z_\ell}{1+z_\ell}
\frac{({(1-\tau z_\ell)^{x_\ell-1}}/{(1-\tau^{-1}z_\ell
)^{x_k-1}})}{(1-\tau z_\ell^2)}(1-z_\ell)^{x_k-x_\ell-1}.
\nonumber
\end{eqnarray*}
Observe that the factors $z_a-z_\ell$ and $1-z_az_\ell$ appearing in
the denominator of
the first line are canceled by matching factors coming out of
$I_{k-1,0}(\vec{x}^{(k)};\vec{z}^{(k)})$. This is crucial, because it
implies that the resulting
integrand has no singularities in $z_\ell$ inside $C_{1,\rho}$ except
possibly at
$z_\ell=1$. On the other hand, since $x_k\geq x_\ell+1$, the
simplification leading to
the second line above implies again that there is no pole at $z_\ell
=1$. We deduce that
the integrand is analytic in $z_\ell$ inside $C_{1,\rho}$, and hence
the integral
vanishes. An analogous argument shows that the residues at $z_k=\tau
^{-1}z_\ell^{-1}$ also
vanish.

Thus, the only important poles are those at $\pm\tau^{-1/2}$. We have
\begin{eqnarray*}
\mathop{\Res}_{z_k=\tau^{-1/2} }I_{k,0}(\vec x;\vec z)& =& I_{k-1,0} \bigl(
\vec x^{(k)}; \vec z^{ (k)}\bigr) \prod
_{a=1}^{k-1} \biggl(\frac{z_a- \tau^{-1/2}}{z_a -
\tau^{1/2}}
\frac{1-\tau^{-1/2}z_a}{1- \tau^{1/2}z_a } \biggr)\\
&&{}\times \frac{
 ( {(1-\tau^{1/2})}/{(1-\tau^{-1/2})}  )^{x_k-1} }{ 2\tau
^{1/2}}
\\
& =& I_{k-1,0} \bigl(\vec x^{(k)}; \vec z^{(k)}\bigr)
(-1)^{x_k-1} \frac{1}2\tau^{({1}/2)x_k-k}.
\end{eqnarray*}
Similarly,
\begin{eqnarray*}
\mathop{\Res}_{z_k=-\tau
^{-1/2}}I_{k,0}(\vec x;\vec z) & =&I_{k-1,0} \bigl(
\vec x^{(k)}; \vec z^{(k)}\bigr) \prod
_{a=1}^{k-1} \biggl(\frac{z_a+ \tau^{-1/2}}{z_a +\tau^{1/2}}
\frac{1+\tau
^{-1/2}z_a}{1+\tau^{1/2}z_a} \biggr)\\
&&{}\times \frac{
 ( {(1+\tau^{1/2})}/{(1+\tau^{-1/2})}  )^{x_k-1} }{
-2\tau^{1/2}}
\\
& =& -I_{k-1,0} \bigl(\vec x^{ (k)}; \vec z^{ (k)}\bigr)
\frac{1}2 \tau^{({1}/2)x_k-k}.
\end{eqnarray*}
If $x_k$ is odd then the two residues cancel each other out. Therefore,
%
%
\begin{eqnarray*}
&&\mathop{\Res}_{z_k=\tau^{-1/2} }I_{k,0}(\vec x;\vec z)+ \mathop{\Res}_{z_k=-\tau^{-1/2}
}I_{k,0}(
\vec x;\vec z) \\
&&\qquad= -I_{k-1,0}\bigl(\vec x^{ (k)};\vec
z^{(k)}\bigr)\mathbf{1}_{x_k\in2\zz_{\geq0}} \tau ^{({1}/2)x_k-k}.
\end{eqnarray*}
Recalling that we have computed the residues on the outside of
$C_{1,\rho}$, which
introduces a minus sign, we get
\begin{eqnarray*}
\tilde u(0,\vec{x})&=&\frac{\tau^{({1}/{2})k(k-1)}}{(2\pi\I
)^{k}}\int_{C_{1,\rho
}^k} \,d\vec z\,
I_{k,0}(\vec z) \\
&= &\tau^{-1}\mathbf{1}_{x_k\in2\zz_{>0}}
\tau^{({1}/2)x_k}\frac{\tau^{({1}/{2})(k-1)(k-2)}}{(2\pi\I
)^{k-1}}\int_{C_{1,\rho}^{k-1}} \,d\vec{z}\,
I_{k-1,0}\bigl(\vec x^{(k)};\vec z\bigr).
\end{eqnarray*}
Equation~\eqref{eq:tildeu0-hf} follows by induction, and this proves (4).
\end{pf*}

\section{Moment formulas}\label{sec:mom}

Recall that Theorem~\ref{Thm:FormulaU} provides a formula for the
expectation of
$\wt{Q}_{x_1}(t)\cdots\wt{Q}_{x_\ell}(t)$, where
$\wt{Q}_{x}(t)=\eta_x(t)\tau^{N_{x-1}(t)}$ and the $x_a$'s have to be
different. To turn
this into a formula for the moments of $\tau^{N_x(t)}$, we will use the
following identity,
first proved as Proposition~3 of \cite{imamSasamASEPduality} (in
\cite{imamSasamASEPduality} the identity was stated only for the
expected value of both
sides, the more general form stated here appears as Lemma~4.17 in \cite{bcs}).
%

\begin{lem}\label{lem:imamSasam}
Let $\eta\in\{0,1\}^\zz$ and write $N_x(\eta)=\sum_{y\leq x}\eta
_y$. Then
%
%
\begin{equation}
\tau^{kN_x(\eta)}=\sum_{\ell=0}^k(-1)^\ell
\pmatrix{k
\cr
\ell} _\tau(\tau ;\tau )_\ell \sum
_{x_1<\cdots<x_\ell\leq x}\eta_{x_1}\tau^{N_{x_1}(\eta
)}\cdots\eta
_{x_\ell}\tau^{N_{x_\ell}(\eta)},\label{eq:imamSasam}
\end{equation}
where the summand for $\ell=0$ should be interpreted as 1.
\end{lem}

Note that this result is not specific to ASEP, which is why we have
introduced the notation
$N_x(\eta)$. For the case of ASEP, and in view of \eqref{eq:Nx}, we
are writing
$N_x(t)=N_x(\eta_t)$. The expected value of the right-hand side of
\eqref{eq:imamSasam} is
explicit in this case (i.e., when we take $\eta$ to be the ASEP
configuration at time $t$, $\eta_t$) thanks to~\eqref{Eq:MomentsQTilde}, and we will turn it into a single multiple
integral it using
arguments similar to those in Section~4 of \cite{bcs}.

\begin{figure}

\includegraphics{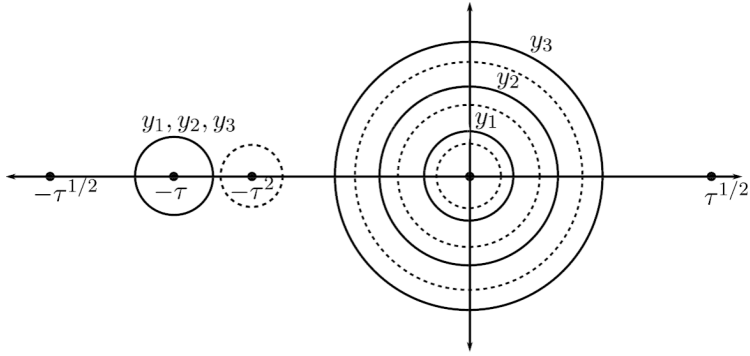}

\caption{Contours appearing in Proposition \protect\ref{prop:NestedPlusSmall}
in the case $k=3$. The
dashed contours correspond to multiplying each of the contours by $\tau
$ and illustrate
the nesting condition described in the proposition.}\vspace*{-5pt}
\label{fig:nested}
\end{figure}

%
\begin{prop} \label{prop:NestedPlusSmall}
For any $k\in\zz_{\geq0}$, we have
%
%
\begin{eqnarray}
\label{Eq:NestedPlusSmall}&& \ee \bigl[\tau^{kN_x(t)} \bigr]
\nonumber
\\[-8pt]
\\[-8pt]
\nonumber
&&\qquad =
\frac{\tau^{({1}/2)k(k-1)}}{(2\pi\I)^k}\int d
\vec{y} \prod_{a<b} \biggl(\frac{y_a-y_b}{y_a-\tau y_b}
\frac{1 - \tau
^{-2}y_ay_b}{1 -
\tau^{-1}y_ay_b} \biggr)\prod_a
\frac{F_{x,t}(y_a)}{y_a},
\end{eqnarray}
where
%
%
\[
\label{Eq:DefF} F_{x,t}(y) = \frac{\tau+y}{\tau-y^2} \biggl(
\frac{1+y}{1+\tau^{-1}y} \biggr)^{x-1} e^{t\hat\ep(y)},
\]
$\hat\ep(y)=\wt\ep(-\tau^{-1}y)$, and the integration contours are
given as follows.
For each $a=1,\ldots,k$, the $y_a$ contour is composed of two
disconnected pieces: a
circle around $-\tau$ with radius small enough so that $-\tau^{1/2}$
is on
its exterior,\vadjust{\goodbreak} and a circle around 0 with radius small enough so that
$\tau^{1/2}$ is on
its exterior. The radii of these circles are chosen so that, in
addition, for all $a<b$
the $y_a$ contour does not include the image under multiplication by
$\tau$ of the $y_b$
contour (see Figure~\ref{fig:nested}).
\end{prop}

%
\begin{pf}
By \eqref{Eq:MomentsQTilde} and Lemma~\ref{lem:imamSasam}, we have
\[
\ee \bigl[\tau^{kN_x(t)} \bigr] =\sum_{\ell=0}^k(-1)^\ell
\pmatrix{k
\cr
\ell} _\tau(\tau;\tau )_\ell G_\ell
\]
with
%
%
\begin{eqnarray}
\label{Eq:Gxkt} G_\ell&=& \frac{\tau^{({1}/{2})\ell(\ell-1)}}{(2\pi\I)^{\ell
}}\int_{C_{1,\rho}^\ell}
\,d\vec{z} \prod_{a<b } \frac{z_a- z_b}{z_a -\tau
z_b }
\frac{1-z_az_b}{1-\tau z_az_b}
\nonumber
\\[-8pt]
\\[-8pt]
\nonumber
&&{}\times \sum_{x_1<\cdots<x_\ell\leq x} \prod
_a \frac{e^{\wt{\ep}(z_a) t}}{ \tau{z_a}^2 -1 } \biggl( \frac{
1-\tau z_a}{ 1-z_a}
\biggr)^{x_a-1}.
\end{eqnarray}
%

For ease of notation, let $\tilde\xi_a=\frac{1-\tau z_a}{1-z_a}$. A
computation shows that
\[
\sum_{x_1<\cdots<x_\ell\leq
x}\prod_{a=1}^\ell
\tilde\xi_a^{x_a-1}=\prod_{a=1}^\ell
\tilde\xi _a^{x-1}\prod_{a=1}^\ell
\frac{1}{\tilde\xi_1\cdots\tilde\xi_a-1}.
\]
Using this in \eqref{Eq:Gxkt}, changing variables $z_a=-\tau^{-1} y_a$
and writing
$\xi_a=\frac{1+y_a}{1+\tau^{-1}y_a}$ we get
%
%
\begin{eqnarray}
\label{Eq:Gxkt2} G_\ell&=&\frac{\tau^{({1}/2)\ell(\ell-1)}}{(2\pi\I)^{\ell }}\int_{C_{-\tau,\tau\rho}^\ell}
\,d\vec{y} \prod_{a< b} \frac{y_a-y_b}{y_a-\tau y_b}\prod
_{a=1}^\ell\frac{1}{\xi_1\cdots
\xi _a-1}
\nonumber
\\[-8pt]
\\[-8pt]
\nonumber
&&{}\times\prod
_{a<b}\frac{1-\tau^{-2}y_ay_b}{1-\tau^{-1}y_ay_b} \prod_{a=1}^\ell
\frac{\xi_a^{x-1}}{\tau-y_a^2}e^{t\hat\ep(y_a)},
\end{eqnarray}
where the new contour $C_{-\tau,\tau\rho}$ is a circle around $-\tau$
with radius
$\tau\rho$ (note that this implies that $-\tau{^{1/2}}$ lies on its
exterior). Now the symmetrization identities appearing in
Lemma~7.2 of \cite{bcs} imply straightforwardly that
\begin{eqnarray*}
&&\sum_{\sigma\in S_\ell} \prod_{a< b}
\frac{y_{\sigma(a)}-y_{\sigma(b)}}{y_{\sigma(a)}-\tau
y_{\sigma(b)}}\prod_a\frac{1}{\xi_{\sigma(1)}\cdots\xi_{\sigma(a)}-1} \\
&&\qquad=
\frac{(-1)^\ell}{(\tau;\tau)_\ell}\prod_a\frac{\tau
+y_a}{y_a}\sum
_{\sigma\in S_\ell} \prod_{a<
b}
\frac{y_{\sigma(a)}-y_{\sigma(b)}}{y_{\sigma(a)}-\tau y_{\sigma(b)}}.
\end{eqnarray*}
Note that, crucially, the last two factors on the right-hand side of
\eqref{Eq:Gxkt2} are
already symmetric, so the above identity can be used to symmetrize the whole
integral, yielding
\[
G_\ell= (-1)^{\ell}\tau^{({1}/2)\ell(\ell-1)-({1}/2)k(k-1)}\frac
{1}{(\tau;\tau)_\ell}
\tilde\nnu_\ell
\]
with
\begin{eqnarray*}
\tilde\nnu_\ell&=&\frac{\tau^{({1}/2)k(k-1)}}{(2\pi\I)^{\ell
}}\int_{C_{-\tau,\tau\rho}^\ell} \,d
\vec{y} \prod_{a<b}\frac
{y_a-y_b}{y_a-\tau
y_b}
\frac{1-\tau^{-2}y_ay_b}{1-\tau^{-1}y_ay_b}\\
&&{}\times \prod_{a=1}^\ell
\frac{e^{t\hat\ep(y_a)}}{\tau-y_a^2}
 \biggl( \frac{1+y_a}{1+\tau^{-1}y_a} \biggr)^{x-1}
\frac{\tau+y_a}{y_a}.
\end{eqnarray*}
Therefore, we have
%
%
\[
\ee \bigl[\tau^{kN_x(t)} \bigr]=\sum_{\ell=0}^k
\pmatrix{k
\cr
\ell } _\tau\tau ^{({1}/2)\ell(\ell-1)-({1}/2)k(k-1)}\tilde
\nnu_\ell.\label
{eq:cauchytype}
\]
We have written things so that we may easily compare with Lemma~4.20 in
\cite{bcs}. Note that
$\tilde\nnu_\ell$ may be rewritten as
\[
\tilde\nnu_\ell=\frac{1}{(2\pi\I)^{\ell}}\int_{C_{-\tau,\tau
\rho}^\ell
}\prod
_{a<b}\frac{y_a-y_b}{y_a-\tau
y_b}s(y_a,y_b)
\prod_af(y_a)\frac{1}{y_a},
\]
where $s(y,y')=(1-\tau^{-2}yy')/(1-\tau^{-1}yy')$ has
no poles in $y$ and $y'$ in a suitable contour encircling $0$ and
$-\tau$, while $f$ is
a function with no poles in a ball around~0 and such that $f(0)=1$.
This is exactly the
structure of $\tilde\nnu_\ell$ in Lemma~4.20 of~\cite{bcs}, and it is
easy to see the extra factor
$\prod_{a<b}s(y_a,y_b)$ in our formula makes no difference in the
argument. Hence,
using their result, we deduce that $\ee [\tau^{kN_x(t)} ]$
has the form claimed in
\eqref{Eq:NestedPlusSmall}.
\end{pf}

As we explained in the \hyperref[sec:intro]{Introduction}, we would
like to manipulate the formula~\eqref{Eq:NestedPlusSmall} given in the last result into one where all
contours coincide.
Doing this involves expanding the nested contours one by one so that
they all end up
coinciding with the largest one. As this multiple contour deformation
is performed, many
poles are crossed.\vadjust{\goodbreak} The associated residues group into clusters, and
this leads to a
formula which is a sum of contour integrals naturally indexed by \emph
{partitions}
$\lambda=(\lambda_1\geq\lambda_{2}\geq\cdots\geq0)$. We will write
$\lambda\vdash k$ if
$\sum_a\lambda_a=k$ and we will denote by $\ell(\lambda)$ the
number of
nonzero elements
of $\lambda$. Additionally, we will write $\lambda
=1^{m_1}2^{m_2}\cdots
$ if $a$ appears
$m_a$ times in $\lambda$, so in this case $\ell(\lambda)=\sum_am_a$ and
$\lambda\vdash\sum_aam_a$.

The contour shift argument referred to above was used in the setting of
Macdonald
processes in \cite{borCor} and later for $q$-TASEP and ASEP in \cite
{bcs}. In the setting
of the delta Bose gas (or Yang's system) with general type root
systems, it goes back to
the work of \cite{heck-op}. Section~7 of \cite{bcpsPlanch-q-TASEP}
contains a detailed
presentation of this argument, and in fact the proposition that follows
is a particular
case of a result proved there.

%
\begin{prop}\label{prop:muk}
\begin{eqnarray}
\label{eq:muk2}&& \ee \bigl[\tau^{kN_x(t)} \bigr]\nonumber\\
&&\qquad= k_\tau!\mathop{
\sum_{
{\lambda\vdash k}}}_{ {\lambda=1^{m_1}2^{m_2}\cdots}}
\frac{(1-\tau)^k}{m_1!m_2!\cdots}
\frac{1}{(2\pi\I)^{\ell(\lambda)}}
\nonumber
\\[-8pt]
\\[-8pt]
\nonumber
&&\qquad\hspace*{71pt}{}\times \int_{\gamma_{-\tau,0}^{\ell(\lambda)}}\,d\vec{w} \det \biggl[
\frac
{-1}{w_a\tau^{\lambda_a}-w_b} \biggr]_{a,b=1}^{\ell(\lambda)}
\\
&&\qquad\hspace*{38pt}\qquad\quad{}\times H\bigl(w_1,\tau w_1,\ldots,
\tau^{\lambda_1-1}w_1,\ldots,w_{\ell(\lambda)},\ldots ,\tau
^{\lambda_{\ell(\lambda)}-1}w_{\ell(\lambda)}\bigr),
\nonumber
\end{eqnarray}
where $\gamma_{-\tau,0}$ is a (positively oriented) contour around
$-\tau$ and $0$,
strictly contained inside the disk of radius $\tau^{1/2}$ and which
does not include any other singularities of the integrand, and
\[
H(y_1,\ldots,y_k) = \prod
_{a<b}\frac{1-\tau^{-2}y_ay_b}{1-\tau^{-1}y_ay_b}\prod_{a}F_{x,t}(y_a).
\]
\end{prop}

\begin{pf}
It is not hard to check that the contours and the integrand which
appear on the right-hand side of \eqref{Eq:NestedPlusSmall} satisfy the
hypotheses of Proposition~7.4 of
\cite{bcpsPlanch-q-TASEP}, and thus
\begin{eqnarray*}
&&\ee \bigl[\tau^{kN_x(t)} \bigr]\\
&&\qquad=\mathop{\sum_{\lambda\vdash
k}}_{ {\lambda=1^{m_1}2^{m_2}\cdots}}
\frac{(-1)^k(1-\tau)^k}{m_1!m_2!\cdots}\frac{1}{(2\pi\I)^{\ell
(\lambda)}} \int_{\gamma_{-\tau,0}^{\ell(\lambda)}}\,d\vec{w} \det
\biggl[\frac
{1}{w_a\tau^{\lambda_a}-w(b)} \biggr]_{a,b=1}^{\ell(\lambda)}
\\
&&\qquad\quad{}\times \prod_{a=1}^{\ell(\lambda)}w_a^{\lambda_a}
\tau^{
({1}/{2})\lambda
_a(\lambda_a-1)}\\
&&\hspace*{28pt}\qquad\quad{}\times E\bigl(w_1,\tau w_1,\ldots,
\tau^{\lambda_1-1}w_1,\ldots,w_{\ell(\lambda)},\ldots ,\tau
^{\lambda_{\ell(\lambda)}-1}w_{\ell(\lambda)}\bigr)
\end{eqnarray*}
with
\begin{eqnarray*}
&&E(y_1,\ldots,y_k)\\
&&\qquad=\sum_{\sigma\in S_k}
\prod_{1\leq b\leq a\leq
k}\frac{y_{\sigma(a)}-\tau y_{\sigma(b)}}{y_{\sigma(a)}-y_{\sigma
(b)}}\prod
_{a<b}\frac{1-\tau^{-2}y_{\sigma(a)}y_{\sigma
(b)}}{1-\tau
^{-1}y_{\sigma(a)}y_{\sigma(b)}} \prod_{a}
\frac{F_{x,t}(y_{\sigma(a)})}{y_{\sigma(a)}}.\vadjust{\goodbreak}
\end{eqnarray*}
Note now that the second and
third products in the definition of $E$ are symmetric under
permutation of the indices in
$\vec y$. On the other hand, by III.(1.4) in \cite{macdonald} the
first double product in the
same identity can be symmetrized as
%
%
\begin{equation}
\label{eq:symmmc} \sum_{\sigma\in S_k}\prod
_{a>b}\frac{y_{\sigma(a)}-\tau
y_{\sigma(b)}}{y_{\sigma(a)}-y_{\sigma(b)}}=(1-\tau)^{-k}(\tau ;\tau
)_k=k_\tau!.
\end{equation}
Hence, $E(y_1,\ldots,y_k)=k_\tau!
\prod_{a<b}\frac{1-\tau^{-2}y_{a}y_{b}}{1-\tau^{-1}y_{a}y_{b}}\prod_{a=1}^{k}\frac{F_{x,t}(y_{a})}{y_{a}}$.
Evaluating $E$ at the point $(y_1,\ldots,y_k)=(w_1,\tau w_1,\ldots
,\tau
^{\lambda_1-1}w_1,\ldots,
w_{\ell(\lambda)},\ldots,\tau^{\lambda_{\ell(\lambda)}-1}w_{\ell
(\lambda)})$ leads,
after some simplifications, to \eqref{eq:muk2}.
\end{pf}

As we will see below, the strings of geometric progressions appearing
in \eqref{eq:muk2}
account for the ratios of $q$-Pochhammer symbols in \eqref
{eq:germans-hf} [see
\eqref{eq:Hsimpl}], which in this case can be thought of as ratios of $q$-Gamma
functions. This is analogous to the strings of arithmetic progressions
which appear in the
case of the delta Bose gas, which give rise to ratios of Gamma
functions (see Section~\ref{sec:app-bose}).

We are now ready for the proof of our main moment formula for $\tau
^{N_x(t)}$ in the
half-flat case.

\begin{pf*}{Proof of Theorem~\ref{thmm:main-hf}}
The formula given in Proposition~\ref{prop:muk} can be rewritten as
%
%
\begin{eqnarray}\label{eq:eetaukNx}
\ee\bigl[\tau^{kN_x(t)}\bigr]&=&k_\tau!\sum
_{\ell=0}^k\mathop{\sum_{m_1,m_2,\cdots}}_{\sum_a
m_a=\ell, \sum_a a m_a=k}
\frac{1}{\ell!}\frac{\ell!}{m_1!m_2!\cdots}\frac{1}{(2\pi\I
)^{\ell}}
\nonumber
\\[-8pt]
\\[-8pt]
\nonumber
&&\hspace*{106pt}{}\times \int
_{\gamma_{-\tau,0}^\ell}\,d\vec w\, I_\ell(\lambda_{m_1,m_2,\ldots};\vec
w),
\end{eqnarray}
where $\lambda_{m_1,m_2,\ldots}$ is
specified by $\lambda_{m_1,m_2,\ldots}=1^{m_1}2^{m_2}\cdots$ and
%
%
\begin{eqnarray}\label{eq:Iell}
&&I_\ell(\lambda;\vec w)\nonumber\\
&&\qquad=\det \biggl[\frac{-1}{w_a\tau^{\lambda_a}-w_b}
\biggr]_{a,b=1}^{\ell(\lambda)}H\bigl(w_1,
\ldots,w_1^{\lambda
_1-1},\ldots ,w_{\ell(\lambda)},\ldots
,w_{\ell(\lambda)}^{\lambda_{\ell(\lambda)}-1}\bigr)\\
&&\qquad\quad{}\times \prod_a(1-
\tau )^{\lambda
_a}.\nonumber
\end{eqnarray}
In the above sum, $m_1,m_2,\ldots$ encodes the partition $\lambda
_{m_1,m_2,\ldots}$ of $k$ of length $\ell$. Observe on the other hand
that, by the symmetry
of the integrand, the right-hand side of \eqref{eq:eetaukNx} is
unchanged if we permute
the $\lambda_a$'s. Thus, we can get rid of the multinomial coefficient\vspace*{1pt}
$\frac{\ell!}{m_1!m_2!\cdots}$ by
replacing the sum over the $m_a$'s by a sum over (unordered)
$n_1,\ldots,n_\ell$ with
the following correspondence: for each $a$, exactly $m_a$ out of the
$n_1,n_2,\ldots,n_\ell$ equal
$a$. This gives
%
%
\begin{equation}\quad
\ee\bigl[\tau^{kN_x(t)}\bigr]=k_\tau!\sum
_{\ell=0}^k \frac{1}{\ell!} \mathop {\sum
_{n_1,\ldots,n_\ell\geq1}}_{ \sum n_a=k} \frac{1}{(2\pi\I)^{\ell}} \int
_{\gamma_{-\tau,0}^\ell}\,d\vec w\, I_\ell \bigl((n_1,
\ldots,n_\ell);\vec w\bigr),\label{eq:mukns}
\end{equation}
where the notation \eqref{eq:Iell} has been extended trivially to unordered
$\ell$-tuples $(n_1,\ldots,n_\ell)$.

What remains is to simplify the integrand. Define
\begin{eqnarray*}
g_1(w)&=&\frac{ (-\tau^{-1}w;\tau )_\infty}{ (\tau
^{-1}w^2;\tau^2 )_\infty} \biggl(\frac{\tau}{\tau+w}
\biggr)^{x-1}e^{(q-p)t({\tau}/{(\tau+w)})},\\
 g_2(w_1,w_2)&=&
\frac{ (\tau^{-1}w_1^2;\tau^2 )_\infty}{ (\tau
^{-3}w_2^2;\tau^2 )_\infty} \frac{ (\tau^{-3}w_2^2;\tau )_\infty}{ (\tau
^{-2}w_1w_2;\tau )_\infty},
\end{eqnarray*}
and write $\vec w\circ{\vec n}=(w_1,\ldots,w_1^{n_1-1},\ldots
,w_{\ell
},\ldots,w_{\ell}^{n_\ell-1})$. We
have
%
%
\begin{eqnarray}\label{eq:Hvecw}
H(\vec w\circ{\vec n})=\wt H(\vec w\circ{\vec n})\prod
_{a=1}^k\prod_{b=0}^{n_a-1}F_{x,t}
\bigl(\tau^bw_a\bigr)
\nonumber
\\[-8pt]
\\[-8pt]
\eqntext{\mbox{with }\displaystyle  \wt H(y_1,
\ldots,y_k)=\prod_{a<b}
\frac{1-\tau
^{-2}y_ay_b}{1-\tau^{-1}y_ay_b}.}
\end{eqnarray}
One checks directly that $F_{x,t}(y)=g_1(y)/g_1(\tau y)$, whence
%
%
\begin{equation}
\prod_{a=1}^k\frac{g_1(w_a)}{g_1(\tau^{n_a}w_a)}= \prod
_{a=1}^k\prod
_{b=0}^{n_a-1}F_{x,t}\bigl(
\tau^bw_a\bigr).\label{eq:ratio-gxt}
\end{equation}
On the other hand, we have
\begin{eqnarray*}
\wt H(\vec w\circ{\vec n})&=&\wt H\bigl(\vec w^{(1)}\circ{\vec
n^{(1)}}\bigr) \prod_{0\leq a_1<a_2<n_1}
\frac{1-\tau^{a_1+a_2-2}w_1^2}{1-\tau
^{a_1+a_2-1}w_1^2}\\
 &&{}\times\prod_{b=2}^k\prod
_{a_1=0}^{n_1-1}\prod
_{a_2=0}^{n_b-1}\frac{1-\tau
^{a_1+a_2-2}w_1w_b}{1-\tau^{a_1+a_2-1}w_1w_b}.
\end{eqnarray*}
The first product on the right-hand side equals
\begin{eqnarray*}
\prod_{a_1=0}^{n_1-2}\frac{1-\tau^{2a_1-1}w_1^2}{1-\tau^{a_1+n_1-2}w_1^2} &=&
\prod_{a_1=0}^{n_1-2}\frac{(\tau^{2a_1-1}w_1^2;\tau^2)_\infty
}{(\tau
^{2a_1+1}w_1^2;\tau^2)_\infty}
\frac{(\tau^{a_1+n_1-1}w_1^2;\tau)_\infty}{(\tau
^{a_1+n_1-2}w_1^2;\tau
)_\infty}\\
& =&g_2\bigl(w_1,\tau^{n_1}w_1
\bigr).
\end{eqnarray*}
One checks similarly that, for fixed $b$, the second product
equals $\mfh(\tau^{-1}w_1,\break\tau^{-1}w_b;  n_1,n_b)$.
We deduce that $\wt H(\vec w\circ{\vec n})=\wt H(\vec w^{(1)}\circ
\vec
n^{(1)})g_2(w_1,\tau^{n_1}w_1) \times\break \prod_{b=2}^k\mfh(\tau^{-1}w_1,  \tau
^{-1}w_b;n_1,n_b)$.
Proceeding inductively to rewrite the right-hand side yields and using
\eqref{eq:Hvecw}
and \eqref{eq:ratio-gxt} yields
%
%
\begin{equation}\qquad
H(\vec w\circ{\vec n})=\prod_a
\frac{g_1(w_a)}{g_1(\tau
^{n_a}w_a)}g_2\bigl(w_a,\tau^{n_a}w_a
\bigr) \prod_{a<b}\mfh\bigl(\tau^{-1}w_a,
\tau^{-1}w_b;n_a,n_b
\bigr).\label{eq:Hsimpl}
\end{equation}
To finish, we note that there is a simplification in the $\tau
$-Pochhammer symbols coming
from the factors $g_1(w_a)/g_1(\tau^{n_a}w_a)$ and $g_2(w_a,\tau^{n_a}w_a)$:
\begin{eqnarray*}
&&\frac{ (-\tau^{-1}w;\tau )_\infty}{ (\tau
^{-1}w^2;\tau^2 )_\infty} \frac{ (\tau^{-1+2n}w^2;\tau^2 )_\infty}{ (-\tau
^{-1+n}w;\tau )_\infty} \frac{ (\tau^{-1}w^2;\tau^2 )_\infty}{ (\tau
^{-3+2n}w^2;\tau^2 )_\infty} \frac{ (\tau^{-3+2n}w^2;\tau )_\infty}{ (\tau
^{-2+n}w^2;\tau )_\infty}
\\
&&\qquad=\frac{ (-\tau^{-1}w;\tau )_\infty}{ (-\tau
^{-1+n}w;\tau )_\infty} \frac{ (\tau^{-1+2n}w^2;\tau^2 )_\infty}{ (\tau
^{-2+n}w^2;\tau )_\infty} \frac{ (\tau^{-3+2n}w^2;\tau )_\infty}{ (\tau
^{-3+2n}w^2;\tau^2 )_\infty} \\
&&\qquad=
\frac{ (-\tau^{-1}w;\tau )_\infty}{ (-\tau
^{-1+n}w;\tau )_\infty} \frac{ (\tau^{-2+2n}w^2;\tau )_\infty}{ (\tau
^{-2+n}w^2;\tau )_\infty}.
\end{eqnarray*}
The right-hand side is exactly $\mfg(w,n)$. Using this in \eqref
{eq:Hsimpl} and
\eqref{eq:Iell}, we deduce that
\begin{eqnarray*}
I_\ell\bigl((n_1,\ldots,n_\ell);\vec w
\bigr)&=&\det \biggl[\frac{-1}{w_a\tau
^{\lambda_a}-w_b} \biggr]_{a,b=1}^\ell
\prod_a\mff\bigl(\tau^{-1}w_a,n_a
\bigr)\mfg\bigl(\tau^{-1}w_a,n_a\bigr)
\\
&&{}\times\prod_{a<b}\mfh\bigl(\tau^{-1}w_a,
\tau^{-1}w_b;n_a,n_b\bigr).
\nonumber
\end{eqnarray*}

Comparing with \eqref{eq:mukns} and \eqref{eq:stpt1} yields the result
after the change
of variables $w_a\mapsto\tau w_a$ (absorbing the Jacobian from the
change of variables
into the determinant).
\end{pf*}

\section{Generating function}\label{sec:genfn}

Since, by definition, $N_x(t)\geq0$, we have $\tau^{N_x(t)}\leq1$ and
thus by
\eqref{eq:qexp} we have for $|\zeta|<1$ that
%
%
\begin{equation}
\label{eq:etau-hf} \ee \bigl[e_\tau\bigl(\zeta\tau^{N_x(t)}\bigr)
\bigr]=\sum_{m\geq0}\frac
{\zeta
^m}{m_\tau!}
\ee^{\mathrm{ h\mbox{-}fl}} \bigl[\tau^{mN_x(t)} \bigr].
\end{equation}
Using \eqref{eq:stpt1} to write the expectation on the right-hand side
explicitly and interchanging the sums in $m$ and $k$ formally leads to
\begin{eqnarray}
\label{eq:etau-hf2} \ee \bigl[e_\tau\bigl(\zeta\tau^{N_x(t)}\bigr)
\bigr]&=&\sum_{k\geq0}\frac
{1}{k!}\sum
_{n_1,\ldots,n_k\geq1}\frac{1}{(2\pi\I)^{k}} \int_{\gamma_{-1,0}^k} \,d\vec
w \det \biggl[\frac{-1}{w_a\tau
^{n_a}-w_b} \biggr]_{a,b=1}^{k}
\nonumber
\\[-8pt]
\\[-8pt]
\nonumber
&&{}\times \prod_a\zeta^{n_a}
\mff(w_a,n_a)\mfg(w_a,n_a)\prod
_{a<b}\mfh (w_a,w_b;n_a,n_b).
\nonumber
\end{eqnarray}
As we will see in the proof of Theorem~\ref{thmm:gen-hf}, the
application of Fubini's
theorem here can be justified, which implies that the above formula
holds as long as
$|\zeta|<1$. In order to analytically extend this identity beyond this
region, we proceed
as in~\cite{borCor} and use a Mellin--Barnes representation for the
sums in $n_a$. The
precise result we will use is the following.

%
\begin{lem}
\label{Lem:MBResult}
Let $g$ be a meromorphic function and $C_{1,2,\ldots}$ a negatively
oriented contour
enclosing all positive integers (e.g., $C_{1,2,\ldots}=\frac{1}{2}+\I
\rr$
oriented with
increasing imaginary part) but no other singularities of $g(\tau^s)$ (in
$s$).\setcounter{footnote}{3}\footnote{Here, $z\longmapsto z^s$ is defined by taking a branch
cut along the
negative real axis.} Then for $\zeta\in\cc\setminus\rr_{\geq0}$ with
$|\zeta|<1$ we
have
\[
\label{Eq:MBResult} \sum_{n=1}^\infty g\bigl(
\tau^n\bigr)\zeta^n=\frac{1}{2\pi\I}\int
_{C_{1,2,\ldots}}\,ds \frac
{\pi}{\sin
(-\pi s)}(-\zeta)^sg\bigl(
\tau^s\bigr),
\]
provided that the left-hand side converges and that there exist closed
contours~$C_k$, $k\in\nn$ enclosing the positive integers from $1$ to
$k$ and such that
the integral of
the integrand on the right-hand side over the symmetric difference of
$C_{1,2,\ldots}$
and $C_k$ goes to zero as $k\to\infty$.
\end{lem}

The statement follows easily from the fact that $\pi/\sin(-\pi s)$ has
a pole at each
$s=k\in\zz$ with residue equal to $(-1)^{k+1}$.

We will also need some precise estimates on $\mfh$, which will be
provided by the lemma
that follows. These estimates will be valid when the relevant variables
lie inside some
carefully chosen contours, which we define next.

\begin{figure}

\includegraphics{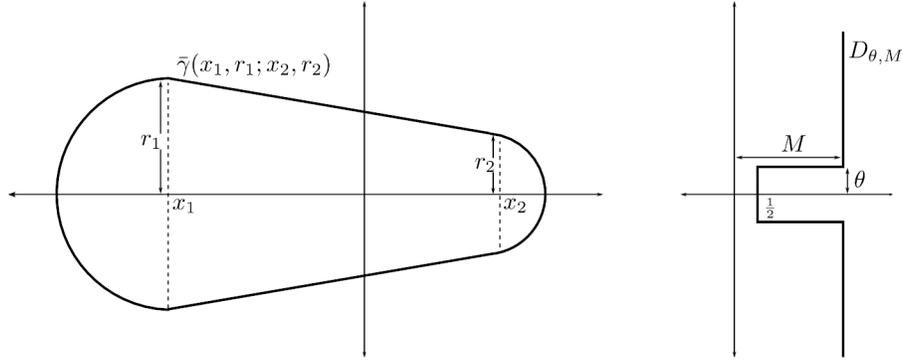}

\caption{Contours in Definition \protect\ref{defn:littleecontours}.}
\label{fig:littleecontour}
\end{figure}

%
\begin{defn}\label{defn:littleecontours}
Let $B(x,r)\subseteq\cc$ denote the ball of radius $r$ centered at
$x$. For $x_1<x_2$
and suitably small $r_1,r_2>0$, we define a positively oriented contour
$\bar\gamma(x_1,r_1;x_2,r_2)$ consisting on the left half of
$\partial
B(x_1,r_1)$, the right
half of $\partial B(x_2,r_2)$, and two lines connecting, respectively, the
top and bottom ends
of the two half circles. Additionally, for $\theta,M>0$ we define a
contour $D_{\theta,M}$
going by straight lines from $M-\I\infty$, to $M-\I\theta$, to
$\frac{1}{2}-\I\theta$, to
$\frac{1}{2}+\I\theta$, to $M+\I\theta$, to $M+\I\infty$. See
Figure~\ref{fig:littleecontour}.
\end{defn}

%

%
\begin{lem}\label{lem:pochinfbd}
Define the function
\[
\mfh_0(z;s_1,s_2)=\frac{ (z;\tau )_\infty (\tau
^{s_1+s_2}z;\tau  )_\infty}{ (\tau^{s_1}z;\tau
)_\infty (\tau^{s_2}z;\tau )_\infty}.
\]
Then there exist constants $C>0$ and $\rho\in(0,\min\{\frac
{1}{2}(\tau
^{-1/2}-1),1\})$ such that, given any
$\delta\in(0,1)$ there are $\theta,M>0$ with the following property:
if $s_1,s_2$ lie to the
right of $D_{\theta,M}$ and $z$ is inside $\bar\gamma(0,\delta
;1,\rho
)$, then
$\llvert \mfh_0(z;s_1,s_2)\rrvert <1+C\delta$.
\end{lem}

\begin{pf}
Fix $\delta_0\in(1,\tau^{-1/2})$
and  $\rho_0\in (0,\min\{\frac12(\tau^{-1/2}-1),1\})$. For fixed $s_1$
and $s_2$,
$\mfh_0(z;s_1,s_2)$ is a meromorphic function of $z$, with poles at
$z=\tau^{-s_1-\ell}$
and $z=\tau^{-s_2-\ell}$ for $\ell\geq0$. Since we are interested
only in
$\Re(s_1)=\Re(s_2)\geq\frac{1}{2}$, all these poles lie outside of
$B(0,\tau
^{-1/2})$, and thus $\mfh_0(z;s_1,s_2)$
is analytic in $z$ inside $\bar\gamma(0,\delta_0;1,\rho_0)$. Now, in
general, if
$D_1,\ldots,D_m$ are bounded domains in $\cc$ and $f$ is a
complex-valued function
defined on $D=D_1\times\cdots\times D_m$ which is analytic in each
variable, then by
the mean value theorem there exists a constant $C>0$ such that for
every $\vec w\in D$
and every $\vec w'\in B(\vec w,\delta)\cap D$ we have
%
%
\begin{equation}
\label{eq:multivar} \bigl\llvert f\bigl(\vec w'\bigr)-f(\vec w)\bigr
\rrvert \leq C\delta.
\end{equation}
We deduce that there is a $C_1>0$ such that if $z,z'$ lie inside
$\bar\gamma(0,\delta_0;1,\rho_0)$ and $|z-z'|<r$, then
%
%
\begin{equation}
\label{eq:C1bd} \bigl|\mfh_0(z;s_1,s_2)\bigr|\leq\bigl|
\mfh_0\bigl(z';s_1,s_2
\bigr)\bigr|+C_1r.
\end{equation}

Now for $x,\alpha_1,\alpha_2,\in[0,1]$ let
\[
g(x;\alpha_1,\alpha_2)=\frac{ (x;\tau )_\infty
(\alpha _1\alpha _2x;\tau )_\infty}{
 (\alpha_1x;\tau )_\infty (\alpha_2x;\tau
)_\infty}.
\]
A computation shows that
$\partial_xg(x;\alpha_1,\alpha_2) |_{x=0}=(\tau
-1)^{-1}(1-\alpha
_1)(1-\alpha_2)$. We deduce that
%
%
\begin{equation}
C_0:=-\sup_{s_1,s_2\in[1/2,\infty)}\partial_xg\bigl(x;
\tau^{s_1},\tau ^{s_2}\bigr) \big|_{x=0}\in(0,
\infty).\label{eq:C0}
\end{equation}
On the other hand, we claim that $g(x;\alpha_1,\alpha_2)$ is concave in
$x\in[0,1]$ for
every fixed $\alpha_1,\alpha_2\in(0,1)$. To see this, write
$g(x;\alpha_1,\alpha_2)=\prod_{\ell\geq0}g_\ell(x;\alpha
_1,\alpha_2)$ with
$g_\ell(x;\alpha_1,\alpha_2)=\frac{(1-\tau^\ell x)(1-\tau^\ell
\alpha
_1\alpha_2
x)}{(1-\tau^\ell\alpha_1x)(1-\tau^\ell\alpha_2 x)}$. Then it is enough
to show that
each $g_\ell$ is positive, decreasing, and concave. The positivity of
$g_\ell$ is clear,
while the decrease and concavity can be checked by computing $\partial
_xg_\ell$ and
$\partial_x^2g_\ell$ (we leave the details to the reader). As a consequence
of this and \eqref{eq:C0}, and since
$\mfh_0(x;s_1,s_2)=g(x;\tau^{s_1},\tau^{s_2})$ and $g(0;s_1,s_2)=1$,
we deduce that
%
%
\begin{equation}
\label{eq:C0bd} \mfh_0(x;s_1,s_2)
\leq1-C_0x
\end{equation}
for all $s_1,s_2\in[\frac{1}{2},\infty)$ and $x\in[0,1]$.

Choose $\rho<\min\{\rho_0,C_0/C_1\}$ and let $r(x)=(1-x)\delta
+x\rho$.
In order to prove
the result it is enough to prove the following statement: there are
$\theta,M>0$
(depending on $\delta$) and $C_2>0$ such that for all $x\in[0,1]$,
$z\in B(x,r(x))$ and
$s_1,s_2$ lying to the right of $D_{\theta,M}$ we have
%
%
\begin{equation}
\label{eq:interp} \bigl|\mfh_0(z;s_1,s_2)\bigr|
\leq1+(C_1+C_2)\delta.
\end{equation}
Assume first that $s_1,s_2\in[\frac{1}{2},\infty)$. Fix $x\in[0,1]$ and
$z\in
B(x,r(x))$. Then by \eqref{eq:C1bd} and \eqref{eq:C0bd}, we have
%
%
\begin{eqnarray}
\label{eq:real-s1s2} \bigl|\mfh_0(z;s_1,s_2)\bigr|&\leq&\bigl|
\mfh_0(x;s_1,s_2)\bigr|+r(x)C_1
\nonumber
\\[-8pt]
\\[-8pt]
\nonumber
&\leq& 1+C_1\delta +C_1(\rho-\delta)x-C_0x<1+C_1
\delta,
\end{eqnarray}
so, in particular, \eqref{eq:interp} holds.

Now we want to extend this to all $s_1,s_2$ lying to the right of
$D_{\theta,M}$. Write
$s_a=\eta_a+\I\theta_a$. There are four cases to consider, depending
on whether or not
$\eta_1$ and $\eta_2$ are larger than $M$. Let us assume first that
$\eta_1,\eta_2\geq
M$. Since $z\in B(0,2)$ (because $\delta,\rho<1$) we have that $\tau
^sz\in
B(0,2\tau^M)\subseteq B(0,\frac{1}{2}\delta)$ for $\Re(s)\geq\frac
{1}{2}$ and
large enough $M$, and
thus $|\tau^{s_1}z-\tau^{\eta_1}z|<\delta$, $|\tau^{s_2}z-\tau
^{\eta
_2}z|<\delta$ and
$|\tau^{s_1+s_2}z-\tau^{\eta_1+\eta_2}z|<\delta$. An argument similar
to the one above,
based on \eqref{eq:multivar}, shows then that there is a constant
$C_2>0$ such that
\[
\bigl\llvert \mfh_0(z;s_1,s_2)-
\mfh_0(z;\eta_1,\eta_2)\bigr\rrvert
<C_2\delta.
\]
Using this together with the bound \eqref{eq:real-s1s2} for $\mfh
_0(z;\eta_1,\eta_2)$
yields \eqref{eq:interp}.

The other three cases are similar. For example, if both $s_1$ and
$s_2$ are in
$[\frac{1}{2},M]\times\I[-\theta,\theta]$ then, for $M$ fixed as
above, we
can choose a
small enough $\theta$ so that $|\tau^{s_1}z-\tau^{\eta_1}z|<\delta$,
$|\tau^{s_2}z-\tau^{\eta_2}z|<\delta$ and
$|\tau^{s_1+s_2}z-\tau^{\eta_1+\eta_2}z|<\delta$, and then the same
argument works. The
mixed case works similarly (although it may yield a different constant).
\end{pf}

\begin{pf*}{Proof of Theorem~\ref{thmm:gen-hf}}
We will prove this result in three steps. The first one will consist
in showing that
\eqref{eq:etau-hf2} holds when $|\zeta|<1$. In the second step we
will apply
the Mellin--Barnes representation given by Lemma~\ref{Lem:MBResult} to turn
\eqref{eq:etau-hf2} into \eqref{eq:gen-hf} for $|\zeta|<1$, $\zeta
\notin\rr_{\geq0}$. Finally we will analytically extend the
resulting formula to all $\zeta\notin\rr_{\geq0}$.

Assume then that $|\zeta|<1$, so that \eqref{eq:etau-hf} holds. Using this
formula together with~\eqref{eq:stpt1} leads to
\[
\ee^{\mathrm{ h\mbox{-}fl}} \bigl[e_\tau\bigl(\zeta\tau^{N_x(t)}\bigr)
\bigr]=\sum_{m\geq
0}\sum_{k=0}^{m}
\frac{1}{k!}\mathop{\sum_{n_1,\ldots,n_k\geq
1}}_{n_1+\cdots+n_k=m}
I_k(\vec n),
\]
where
%
%
\begin{eqnarray*}
I_k(\vec n)&=&\frac{1}{(2\pi\I)^{k}} \int_{\gamma_{-1,0}^k} \,d\vec
w \det \biggl[\frac{-1}{w_a\tau
^{n_a}-w_b} \biggr]_{a,b=1}^{k}
\\
&&{}\times \prod
_a\zeta^{n_a}\mff(w_a;n_a)
\mfg(w_a;n_a)\prod_{a<b}\mfh
(w_a,w_b;n_a,n_b).
\end{eqnarray*}
Interchanging the sums in $k$ and $m$ leads to
%
%
\begin{eqnarray}\label{eq:gen-hf-pf}
\ee^{\mathrm{ h\mbox{-}fl}} \bigl[e_\tau\bigl(\zeta\tau^{N_x(t)}\bigr)
\bigr] &= &\sum_{k\geq
0}\sum_{
m\geq k}
\frac{1}{k!}\mathop{\sum_{n_1,\ldots,n_k\geq
1}}_{n_1+\cdots+n_k=m}
I_k(\vec n)
\nonumber
\\[-8pt]
\\[-8pt]
\nonumber
&= &\sum_{k\geq0}
\frac{1}{k!}\sum_{n_1,\ldots,n_k\geq
1}I_k(\vec
n).
\end{eqnarray}
In order to justify the application of Fubini's theorem, it is enough
to verify that the
sum $\sum_{k\geq0}\sum_{m\geq k} |\frac{1}{k!}\sum_{n_1+\cdots+n_k=m}
I_k(\vec n) |$
is finite, which by the triangle inequality, will follow if we verify that
%
%
\begin{equation}
\label{eq:fubini-check} \sum_{k\geq0}\frac{1}{k!}\sum
_{n_1,\ldots,n_k\geq1} \bigl|I_k(\vec n)\bigr|<\infty.
\end{equation}
The main difficulty we face at this point is the fact that the absolute
value of
$\mfh(w_a,w_b;n_a,n_b)$ is in general not bounded by 1, which in
principle introduces a
factor of order $c^{k^2}$ into our sum for some $c>1$. To deal with
this issue, we will
have to choose the contour $\gamma_{-1,0}$ carefully, and moreover let
it depend on
$k$. Note, however, that this choice is made at this point only in
order to obtain a
suitable estimate, and does not fix the contour in the statement of
the theorem.

Now fix $\rho>0$ and $C>0$ as in Lemma~\ref{lem:pochinfbd} and, for
fixed $k$, let
$\delta_k=C^{-1}(2^{1/k}-1)$ and choose $\theta_k,M_k>0$ as in Lemma~\ref{lem:pochinfbd} for $\delta=\delta_k$. Furthermore, let $\delta
_k',\rho'>0$, $\theta_k'<\theta_k$ and $M_k'>M_k$,
and write $\bar\gamma_k=\bar\gamma(-1,\rho';0,\delta_k')$ and
$\bar
D_k=D_{\theta_k',M_k'}$ ($\bar D_k$ will be used in the second step).
Note that
$\bar\gamma_k$ is star-shaped with respect to the origin (i.e., any
ray emanating from
the origin intersects the contour in one and only one point). This
implies, in
particular, that the denominator inside the determinant appearing in
$I(\vec n)$ never
vanishes. On the other hand, by choosing $\delta_k'$ and $\rho'$ to be
suitably small we
may assume that $\bar\gamma_k$ is contained inside $B(0,\tau^{-1/2})$,
in which case it
is easy to check that there are not singularities of $\mfh$ inside.
Therefore, our choice
of $\bar\gamma_k$ satisfies the requirements of Theorem~\ref{thmm:main-hf}.

Having made this choice of contour, we claim that we can choose an
$\eta>0$ such that if
$\delta_k'=\eta\delta_k$ and $\rho'$ is small enough then whenever
$w_a,w_b\in\bar\gamma_k$ we have that $w_aw_b$ is contained inside
$\bar\gamma(0,\delta_k;1,\rho)$. To see this, observe that
$\{ww' \dvtx w,w'\in[-1,0]\}=[0,1]$ and, therefore, given any open
neighborhood $U$ of $[0,1]$
we can find an open neighborhood $V$ of $[-1,0]$ such that $\{ww'
\dvtx w,w'\in V\}$ is
contained inside $U$. Our claim follows easily from this because given
any such
neighborhood $V$ we can choose $\delta_k'$ and $\rho$ small enough so that
$\bar\gamma_k$ is contained inside $V$.

Making these choices, and thanks to our earlier choices of parameters
and using Lemma~\ref{lem:pochinfbd}, we get
%
%
\begin{equation}
\bigl\llvert \mfh(w_a,w_b;n_a,n_b)
\bigr\rrvert =\bigl\llvert \mfh _0(w_aw_b;n_a,n_b)
\bigr\rrvert \leq2^{1/k}\label{eq:mfhbd}
\end{equation}
for $w_a,w_b\in\bar\gamma_k$ and $n_a,n_b\in\zz_{\geq1}$ (since
in this
case $n_a$ and
$n_b$ trivially lie to the right of $\bar D_k$). On the other hand,
the only
singularity of $\mff(w_a;n_a)$ occurs at $w_a=-1$, and since $\bar
\gamma_k$ stays at
distance at least $\rho'$ from $-1$, this factor is uniformly bounded
along the contour,
say by some constant $c_1>0$ (independently of $k$). A~similar
argument shows that $|\mfg(w_a;n_a)|$ is uniformly
bounded (say by $c_1$ again), and we deduce that
%
%
\begin{eqnarray}
\label{eq:Iknbd} %
\bigl\llvert I_k(\vec n)\bigr
\rrvert &\leq&\frac{c_1^{2k}2^{
({1}/{2})(k-1)}}{(2\pi
)^k}\int_{\bar\gamma_k^k} \,d\vec w \prod
_a|\zeta|^{n_a}\biggl\llvert \det \biggl[
\frac{-1}{w_a\tau
^{n_a}-w_b} \biggr]_{a,b=1}^{k}\biggr\rrvert
\nonumber\\
&\leq& c_2^k|\zeta|^{\sum_an_a}\int
_{\bar\gamma_k^k} \,d\vec w \prod_a
\frac{1}{|w_a|}\biggl\llvert \det \biggl[\frac{-w_a}{w_a\tau
^{n_a}-w_b}
\biggr]_{a,b=1}^{k}\biggr\rrvert
\\
&\leq& c_2^kk^{k/2}|\zeta|^{\sum_an_a}\int
_{\bar\gamma_k^k} \,d\vec w \prod_a
\frac{1}{|w_a|}\sup_{a,b=1,\ldots,k}\biggl\llvert \frac
{w_a}{w_a\tau
^{n_a}-w_b}
\biggr\rrvert ^k\nonumber
\end{eqnarray}
for some $c_2>0$, where in the last inequality we used Hadamard's
bound. The supremum is
clearly bounded by some constant $c_3>0$, uniformly in $w_a$, $w_b$
and $n_a$. On the
other hand, it is not hard to check that
\[
\int_{\bar\gamma_k} \,dw_a \frac{1}{|w_a|}\leq
c_4\bigl| \log\bigl(\delta_k'
\bigr)\bigr|=c_4\bigl| \log(\eta\delta_k)\bigr|\leq
c_4'\log(k)
\]
for some $c_4,c_4'>0$ by our choice of $\delta_k$ and $\delta_k'$. We
deduce that
\[
\bigl|I_k(\vec n)\bigr|\leq c^k\bigl(k^{1/2}\log(k)
\bigr)^k|\zeta|^{\sum_an_a}
\]
for some $c>0$ and thus, since we are taking
$|\zeta|<1$, \eqref{eq:fubini-check} holds. Therefore, \eqref
{eq:gen-hf-pf} holds
for $|\zeta|<1$.

As we mentioned at the beginning of the proof, the next step is to
apply the
Mellin--Barnes representation to \eqref{eq:gen-hf-pf}. The idea is to
focus on the $k$th
term of the sum on the right-hand side of \eqref{eq:gen-hf-pf} for
some fixed $k$, and
then apply Lemma~\ref{Lem:MBResult} one by one to each of the sums in
$n_1,\ldots,n_k$
with the contour $C_{1,2,\ldots}$ taken as $\bar D_k=D_{\theta
_k',M_k'}$ [and
$\gamma_{-1,0}$ as $\bar\gamma_k=\bar\gamma(-1,\rho';0,\delta'_k)$],
which would prove
the identity
%
%
\begin{eqnarray}
\label{eq:gen-hf-mod} \ee^{\mathrm{ h\mbox{-}fl}} \bigl[e_\tau\bigl(\zeta
\tau^{N_x(t)}\bigr) \bigr]& =&\sum_{k=0}^{\infty}
\frac{1}{k!}\frac{1}{(2\pi\I)^{2k}}\int_{\bar
D_{k}^k}\,d\vec s \int
_{\bar\gamma_k^k}\,d\vec w\det \biggl[\frac{-1}{w_a\tau
^{s_a}-w_b}
\biggr]_{a,b=1}^{k}
\nonumber
\\[-8pt]
\\[-8pt]
\nonumber
&&{}\times\prod_a(-\zeta)^{s_a}
\mff(w_a;s_a)\mfg(w_a;s_a)\prod
_{a<b}\mfh (w_a,w_b;s_a,s_b)
\end{eqnarray}
for $\zeta\notin\rr_{\geq0}$ with $|\zeta|<(1-\tau)^{-1}$. To
this end,
we need to
verify that the conditions of the lemma are satisfied. Note that, in
view of the
preceding argument, we are free to choose $\theta_k'$ and $M_k'$ to be
respectively even
smaller and even larger than in our original choice. We start by
observing that
$w_a\tau^{s_a}-w_b$ never vanishes for $s_a$ along this contour. To
see this, note first
that $M_k'$ can be chosen to be sufficiently large so that if $\bar
\gamma_k$ is scaled
by $\tau^{M_k'}$ then any rotation of the resulting contour is
contained inside
$\bar\gamma_k$, which shows that $w_a\tau^{s_a}-w_b\neq0$ for $s_a$ with
$\Re(s_a)\geq M_k'$. On the other hand, since $\bar\gamma_k$ is star-shaped,
$w_a\tau^{s_a}-w_b\neq0$ for $s_a\in[\frac{1}{2},\infty)$, and
thus the same
holds in the strip
$[\frac{1}{2},M_k']\times\I[-\theta_k',\theta_k']$ if $\theta_k'$
is small
enough. This shows
that there are no singularities of the determinant in the integrand in
\eqref{eq:gen-hf-mod} for $s_a$ lying to the right of $\bar D_k$. The
singularities of
the remaining factors are all avoided in this region for similar reasons.

What is left to check is that there are closed contours $C_{k,m}$
enclosing $1,\ldots,m$
[and contained in $\{s \dvtx\Re(s)\geq\frac{1}{2}\}$] such that the
integral on
the symmetric
difference of $\bar D_k$ and $C_{k,m}$ goes to 0 as $m\to\infty$. We
choose $C_{k,m}$ to
be union of the piece of $\bar D_k$ lying inside $B(0,m+\frac{1}{2})$ and
the arc on the
boundary of this ball lying to the right of $\bar D_k$. But this is
actually not hard to
see. We have already checked that $\mff(w_a;s_a)$, $\mfg(w_a;s_a)$,
$\mfh(w_a,w_b;s_a,s_b)$ and the determinant have no singularities for
$s_a,s_b$ lying to
the right of $\bar D_k$, and since these factors depend on $s_a,s_b$
only through
$\tau^{s_a},\tau^{s_b}$, which live in a compact set, they are bounded
uniformly. The
necessary decay is going to come from the product $|\pi/\sin(\pi
s_a)||\zeta^{s_a}|$. In
fact, as $|\Im(s_a)|\to\infty$ with $\Re(s_a)=\frac{1}{2}$ we have that
$|\pi
/\sin(\pi s_a)|$
decays exponentially while $|\zeta^{s_a}|$ stays bounded. The same
exponential decay
applies in the circular part of $C_{k,m}$ restricted to $| \arg
(s_a)|>\frac{\pi}{4}$
[since here $|\Im(s_a)|\to\infty$ as before]. Finally, note that on
the circular piece
of $C_{k,m}$ with $| \arg(s_a)|>\frac{\pi}{4}$ we have that $s_a$
stays bounded away from
all integers, so that $|\pi/\sin(\pi s_a)|$ is uniformly bounded,
while $\Re(s_a)\to\infty$, so that
$|\zeta^{s_a}|$ decays exponentially. Putting these facts together
shows that the
integrand has the right decay, and gives \eqref{eq:gen-hf-mod}.

Our third step is to analytically extend \eqref{eq:gen-hf-mod} to all
$\zeta\notin\rr_{\geq0}$, for which we need to show that both sides
are analytic in
$\zeta$ in that region. Observe first that the left-hand side is given by
\[
\ee^{\mathrm{ h\mbox{-}fl}} \bigl[e_\tau\bigl(\zeta\tau^{N_x(t)}\bigr)
\bigr] =\sum_{n\geq0}\frac{\pp^{\mathrm{ h\mbox{-}fl}} (N_x(t)=n
)}{ ((1-\tau )\zeta\tau^n;\tau )_\infty}.
\]
For each $\zeta\notin\{(1-\tau)^{-1}\tau^{-m}\}_{m\in\zz_{\geq0}}$,
this series is
uniformly convergent on a neighborhood of $\zeta$, and thus the
left-hand side is
analytic for $\zeta\notin\rr_{\geq0}$.

Turning to the right-hand side of \eqref{eq:gen-hf-mod}, observe that
each summand in the series is clearly analytic in
$\zeta\notin\rr_{\geq0}$. We will use now the fact that the limit
of a uniformly
absolutely convergent series of analytic functions is analytic to show
that the right-hand side of
\eqref{eq:gen-hf-mod} is analytic in $\zeta$ in any fixed neighborhood
which avoids
$\rr_{\geq0}$. Consider the $k$th term of our series and recall that
we have chosen $\delta_k'$ and $\rho'$ so that
$w_aw_b$ is inside $\bar\gamma(0,\delta_k;1,\rho)$ for $w_a,w_b\in
\bar
\gamma_k$, while
on the other hand $\theta_k'<\theta_k$ and $M_k'>M_k$. As a
consequence, and thanks to
Lemma~\ref{lem:pochinfbd} and our choice of parameters, we deduce as
in \eqref{eq:mfhbd} that
$|\mfh(w_a,w_b;s_a,s_b)|\leq2^{1/k}$ for $w_a,w_b\in\bar\gamma_k$ and
$s_a,s_b\in\bar
D_k$. As in the previous step, we have that $\mff(w_a;s_a)$,
$\mfg(w_a;s_a)$, $\mfh(w_a,w_b;s_a,s_b)$ are uniformly bounded and
proceeding as in
\eqref{eq:Iknbd} we deduce that the $k$th term of the series on the
right-hand side of
\eqref{eq:gen-hf-mod} is bounded in absolute value by
\begin{eqnarray*}
&&\frac{c_1^k}{k!}\frac{1}{(2\pi\I)^{2k}}\int_{\bar D_k^k}\,d\vec s \int
_{\bar\gamma_k^k}\,d\vec w \prod_a\biggl
\llvert \frac{\pi}{\sin
(\pi
s_a)}\biggr\rrvert \frac{|\zeta^{s_a}|}{|w_a|}\sup
_{a,b=1,\ldots,k}\biggl\llvert \frac
{w_a}{w_a\tau^{s_a}-w_b}\biggr\rrvert
^{k}
\\
&&\qquad\leq\frac{c_2^k(k^{1/2}\log(k))^k}{k!}\frac{1}{(2\pi\I)^{k}}\int_{D_k^k}\,d\vec s
\prod_a\biggl\llvert \frac{\pi}{\sin(\pi s_a)}\biggr
\rrvert \bigl|\zeta^{s_a}\bigr| \leq\frac{c_3^k(k^{1/2}\log(k))^k}{k!}
\end{eqnarray*}
for some constants $c_1,c_2,c_3>0$ which are uniform in $\zeta$ in a
compact subset of
$\cc$ [here we have used again the fact that $|\pi/\sin(\pi s_a)|$
decays exponentially
as $\Im(s_a)\to\infty$]. This shows that the right-hand side of
\eqref
{eq:gen-hf-mod} is
absolutely summable, uniformly in $\zeta$ on a fixed neighborhood away from
$\rr_{\geq0}$ as required, and thus finishes the analytic extension of
\eqref{eq:gen-hf-mod} to all $\zeta\notin\rr_{\geq0}$.

At this point, we have proved \eqref{eq:gen-hf-mod}. We may now deform the
contours $\bar D_k$ and $\bar\gamma_k$ in each of the summands to the
contours $\frac{1}2+\I\rr$ and
$\gamma_{-1,0}$ by appealing to Cauchy's theorem, thus finishing the proof.
\end{pf*}

\section{Formulas for the KPZ/stochastic heat equation}\label{sec:app-bose}

The one-dimensional Kardar--Parisi--Zhang (KPZ) ``equation'' is given by
\[
\partial_t h =\tfrac{1}2\partial_x^2
h -\tfrac{1}2 \bigl[(\partial_x h)^2 -\infty
\bigr] +\xi,
\]
where $\xi$ is a space--time white noise. This SPDE is ill-posed as
written but, at least
on the torus, it can be made sense of by a renormalization procedure
introduced by
Hairer in \cite{hairer,hairerReg}. His solutions coincide with the
Cole--Hopf solution
(which is known to be the physically relevant solution; see, e.g., the review
\cite{quastelCDM}) obtained by setting
%
%
\begin{equation}
h(t,x)=-\log Z(t,x),\label{eq:CH}
\end{equation}
where $Z$ is the unique solution to the (well-posed) stochastic heat
equation (SHE)
%
%
\begin{equation}
\label{SHE} \partial_t Z = \tfrac{1}2
\partial_x^2 Z +\xi Z.
\end{equation}
We will now give a contour integral ansatz for the moments of $Z$ with
the ``tilted''
half-flat initial data defined by $Z(0,x)=e^{-\theta x}\mathbf
{1}_{x\geq0}$.

To be more precise, we will provide a solution for the \emph{delta Bose
gas} with this initial data, which we
interpret as the solution $v(t,x)$ to the following system of
equations, where we write
$W_k=\{\vec x\in\rr^k \dvtx x_1<x_2<\cdots<x_k\}$ (see \cite
{borCor} for
more details):
\begin{longlist}[(3)]
\item[(1)] For $\vec x\in W_k$,
\[
\partial_tv(t,\vec x)=\tfrac{1}2\Delta v(t,\vec x),
\]
where the Laplacian acts on $\vec{x}$.
\item[(2)] For $\vec x$ on the boundary of $W_k$, with $x_a=x_{a+1}$,
\[
(\partial_{x_a}-\partial_{x_{a+1}}-1)v(t,\vec x)=0.
\]
\item[(3)] For $\vec x\in W_k$,
\[
\lim_{t\to0}v(t,\vec x)=v_0(\vec x).
\]
\end{longlist}
In the (tilted) half-flat case, we take $v_0(\vec x)=\prod_ae^{-\theta
x_a}\mathbf{1}_{x_a\geq0}$.

It is widely accepted in the physics literature that, if $Z(t,x)$ is a
solution of the
SHE, then $v(t;\vec x)=\ee[Z(t,x_1)\cdots Z(t,x_k)]$ is a solution of
the delta Bose gas. This
fact is proved in \cite{mqrScaling}, where it is also shown that there
is at most one
solution. Therefore, our formulas below for the solution of the delta
Bose gas are indeed
identifying the $\ee[Z(t,x_1)\cdots Z(t,x_k)]$. In any case, in the
last result of this
section (Proposition~\ref{prop:DBGHF-shifted-2}) we will state a
formula both for the
delta Bose gas and for the moments of the solution of the SHE, with a
proof for the second
part which is independent of this correspondence.

Given $\alpha\in\rr^k$, we will write
$\vec\alpha+(\I\rr)^k=(\alpha_1+\I\rr)\times\cdots\times
(\alpha_k+\I\rr)$.

%
\begin{prop}\label{prop:deltaBoseFormulaU}
The delta Bose gas with tilted half-flat initial condition given by
$v_0(\vec x)=\prod_ae^{-\theta x_a}\mathbf{1}_{x_a\geq0}$,
$\theta\geq0$, is solved by
%
\begin{eqnarray}
\label{Eq:DBGHF} v(t,\vec{x})&=& \frac{1}{(2\pi\I)^k} \int_{\vec\alpha+(\I\rr)^k}\,d
\vec z \prod_{a<b} \biggl(\frac{z_a-z_b}{z_a- z_b-1}
\frac{z_a+z_b-1}{z_a+z_b} \biggr)
\nonumber
\\[-8pt]
\\[-8pt]
\nonumber
&&{}\times\prod_{a=1}^k
\frac{1}{z_a}e^{({t}/{2})\sum_{a=1}^k
(z_a-\theta)^2+\sum_{a=1}^k(z_a-\theta)x_a},
\end{eqnarray}
where $\alpha_1>\alpha_2+1>\cdots>\alpha_k+k-1>k-1$ and $x_1<\cdots<x_k$.
\end{prop}

\begin{pf} We only verify that (3) is satisfied, the rest follows as
in the case of
$\delta_0$ initial condition \cite{borCor}. We need to show that
\[
\lim_{t\to0}v(t,\vec x)=\prod_{a-1}^ke^{-\theta x_a}
\mathbf {1}_{x_a\geq0}.
\]
We will denote the integrand by $I_k(z_1,\ldots,z_k)$. Assume first
that $x_1<0$. Thanks
to the factor $e^{({1}/2)t(z_1-\alpha)^2}$ we may move the $z_1$
contour to
$\alpha_1+R+\I\rr$, $R>0$. Note that we do not cross any poles.
Changing variables
$z_1\mapsto z_1+R$ gives
\[
v(t,\vec{x}) = \frac{1}{(2\pi\I)^k} \int_{\vec\alpha+\I\rr}\,d\vec{z}
I_k(z_1+R,z_2,\ldots,z_k).
\]
Now me may compute the limit $t\to0$, which removes the quadratic term
in the
exponential. The resulting integrand in $\lim_{t\to0}v(t,\vec x)$
contains a factor
$\frac{1}{z_1+R}e^{x_1(z_1+R)}$, and since $x_1<0$, we may take $R\to
\infty$ to deduce
without difficulty that the integral vanishes in this case.

So we assume now that $x_1\geq0$ (and so $x_a\geq0$ for all
$a=1,\ldots
,k$). Our goal is
to move the $z_k$ contour to $-M+\I\rr$ (with $M>\alpha_1$). We may do
this thanks to
the Gaussian factor as before. Observe that the poles for $z_k$ on $\{
-M\leq\Re(z_k)\leq
\alpha_k\}$ are $0$ and $-z_a$ for $a<k$. We begin with the second
type of pole. We
have, for $\ell<k$,
\begin{eqnarray*}
&&\mathop{\Res}_{z_k=-z_\ell}I_k(z_1,\ldots,z_k)\\
&&\qquad=
\int \,d\vec z\, I_{k-1}(z_1,\ldots,z_{k-1})
\frac{2z_\ell e^{-(z_\ell+\alpha)x_k+({1}/2)t(z_k-\alpha
)^2}}{2z_\ell
-1}\frac{1}{-z_\ell}
\\
&&\qquad\quad{}\times\mathop{\prod_{a=1}}_{a\ne\ell}^{k-1}
\frac{z_a+z_\ell
}{z_a+z_\ell-1} \frac{z_a-z_\ell-1}{z_a-z_\ell}
\\
&&\qquad=\int \,d\vec z\, I_{k-2}(z_1,\ldots,z_{\ell-1},z_{\ell+1},
\ldots,z_{k-1})\\
&&\qquad\quad{}\times\frac
{-2e^{-z_\ell
(x_k-x_\ell)-\alpha(x_\ell+x_k)
+({1}/2)t(z_\ell-\alpha)^2+({1}/2)t(z_k-\alpha)^2}}{z_\ell(2z_\ell-1)}
\\
&&\qquad\quad{}\times\prod_{b=\ell+1}^{k-1}\frac{1+z_\ell-z_b}{z_\ell-z_b-1}.
\end{eqnarray*}
Observe that, due to the cancellation leading to the second line, the
$z_\ell$ integral
has no poles on $\{\Re(z_\ell)>\alpha_\ell\}$. As before we may freely
move the $z_\ell$ contour to
$\alpha_\ell+R+\I\rr$, $R>0$. Changing variables $z_\ell\mapsto
z_\ell
+R$ and taking
$t\to0$ yields an integral over the original $z_1,\ldots,z_{k-1}$
contours and
containing a factor $e^{-(z_\ell+R)(x_k-x_\ell)-\alpha(x_\ell+x_k)}$
and no quadratic
term in the exponent. Since $x_k>x_\ell$, taking $R\to\infty$ shows
that this term vanishes.

We still need to compute the pole at $z_k=0$, but let us first observe
that the $z_k$
integral over the new contour $-M+\I\rr$ also vanishes after taking
the limit
$t\to0$. In fact, proceeding as above, now changing variables $z_k\to
z_k-M$, the resulting $k$-fold integral equals
\[
v(t,\vec{x}) = \frac{1}{(2\pi\I)^k} \int_{\vec\alpha+\I\rr}\,d\vec{z}
\,I_k(z_1,z_2,\ldots,z_k-M).
\]
In the limit $t\to0$, the integrand contains a factor of the form
$e^{x_k(z_k-M)}$, and
since we are assuming $x_k>0$ we may take $M\to\infty$ to deduce that
the whole integral
goes to 0.

So the only term left in the limit $t\to0$ is the one corresponding to
the pole at
$z_k=0$. We have
\begin{eqnarray*}
&&\mathop{\Res}_{z_k=0}I_k(z_1,\ldots,z_k)\\
&&\qquad=
\int_{\alpha_j+\I\rr} \,d\vec z\, I_{k-1}(z_1,
\ldots,z_{k-1})e^{({1}/2)t\alpha^2-\alpha
x_k}\prod_{a=1}^{k-1}
\biggl(\frac{z_a}{z_a-1} \frac
{z_a-1}{z_a} \biggr).
\end{eqnarray*}
The last product is obviously 1, so we have proved that
\begin{eqnarray*}
\lim_{t\to0}v(t,\vec x)&=&\lim_{t\to0}
\mathbf{1}_{x_k\geq0}\int_{\alpha
_j+\I\rr} \,d\vec z
\,I_{k-1}(z_1,\ldots,z_{k-1})e^{({1}/2)t\alpha^2-\alpha x_k}\\
& =&
\mathbf{1}_{x_k\geq0}e^{-\alpha x_k}\lim_{t\to0}v
\bigl(t,(x_1,\ldots,x_{k-1})\bigr).
\end{eqnarray*}
The result follows by induction.
\end{pf}

Observe that, as should be expected, multiplying \eqref{Eq:DBGHF} by
$\theta^k$ and
letting $\theta\to\infty$ yields (after shifting contours by $\theta$
and changing
variables $z_a\mapsto z_a+\theta$) the solution of \cite{borCor} for
the delta Bose gas
with \emph{narrow wedge} initial condition [which corresponds to
$Z(0,x)=\delta_0(x)$ at
the level of the SHE], given by
%
%
\begin{equation}
\label{Eq:DBGHFlim}\qquad v_0(t,\vec x)= \frac{1}{(2\pi\I)^k}\int
_{\vec\alpha+(\I\rr
)^k}\,d\vec z \prod_{a<b}
\frac{z_a-z_b}{z_a- z_b-1} \prod_{a=1}^ke^{({t}/{2})\sum_{a=1}^k z_a^2+\sum_{a=1}^kz_ax_a}
\end{equation}
for $x_1<\cdots<x_k$.

When $\theta=0$, \eqref{Eq:DBGHF} gives the solution for the half-flat
initial condition $Z(0,x)=\mathbf{1}_{x\geq0}$, which can also be
obtained by
taking the weakly asymmetric limit of
\eqref{Eq:MomentsQTilde} (see the proof of Proposition~\ref
{prop:DBGHF-shifted-2} for a
similar computation).

By linearity of \eqref{SHE}, we have that, if $Z(0,y;t,x)$ is the
solution to the SHE with
initial data $Z(0,y;0,x)=\delta_y(x)$, then $Z(t,x)=\int_{-\infty
}^\infty \,dy Z(0,y,t,x)f(y)$ solves the SHE
with initial condition $Z(0,x)=f(x)$, and hence
\[
\ee_f \bigl[Z(t,x_1)\cdots Z(t,x_k)
\bigr] = \int_{\rr^k}\,d\vec y \,\ee \biggl[\prod
_a Z(0,y_a;t,x_a) \biggr] \prod
_a f(y_a)
\]
(with the subscript in $\ee_f$ denoting the initial condition for the
SHE). Although we do
not have a formula for the integrand on the right-hand side in general
note that, by
statistical time reversal invariance, we do have
\[
\ee \biggl[\prod_a Z(0,y_a;t,x_a)
\biggr] = \ee \biggl[\prod_a Z(0,x_a;t,y_a)
\biggr].
\]
Now if all the $x_a$'s are the same, we can use the spatial statistical
invariance and
symmetry to see that $\ee [\prod_a Z(0,x_a;t,y_a) ]=
\ee [\prod_aZ(t,x-y_a) ]$. Finally, changing variables and then
restricting to the
\emph{Weyl chamber} $W_k=\{\vec{x}\in\rr^k \dvtx x_1<\cdots<x_k\}
$, we obtain
\[
\label{eq:db-moments} \ee_f \bigl[Z(t,x)^k \bigr]=k!\int
_{W_k}\,d\vec y \,\ee \biggl[\prod_aZ(t,y_k-x)
\biggr]\prod_af(x-y_a).
\]
Specializing to the tilted half-flat initial condition given by
$f=f_\theta$ with
$f_\theta(x)=e^{-\theta x}\mathbf{1}_{x\geq0}$, and in view of the relation
between the delta
Bose gas and the moments of the SHE discussed above, this suggests an
alternative route
for obtaining a formula for $v(t;x,\ldots,x)$ in this case, namely
%
%
\begin{equation}
\label{eq:deltaBose-lin} v(t;x,\ldots,x)=k!\int_{W_k}\,d\vec y\,
v_0(t;\vec y)\prod_ae^{-\theta
(x-y_a)}
\mathbf{1}_{y_a\leq x}
\end{equation}
with $v_0$ as in \eqref{Eq:DBGHFlim}. Although this identity can be
justified directly
from the linearity of the delta Bose gas itself, it is not at all clear
at a first look
that this alternative computation would lead to the same formula as the
one in Proposition~\ref{prop:deltaBoseFormulaU}.

To see directly why the above formula holds, we start by using the
explicit formula for $v_0(t;\vec
y)$ and computing the $y_a$ integrals over $W_k$, which yield
\[
\frac{k!}{(2\pi\I)^k}\int_{\vec\alpha+(\I\rr)^k}\,d\vec z \prod
_{a<b}\frac{z_a-z_b}{z_a- z_b-1} \prod_{a=1}^k
\frac{1}{z_1+\cdots+z_a}e^{({t}/{2})\sum_{a=1}^k
z_a^2+\sum
_{a=1}^kz_ax}.
\]
Now deform the $z_a$ contours one by one so that they all coincide with
the leftmost
one. The answer is obtained by an argument analogous to the proof of Proposition~\ref{prop:muk}, and is given by
\begin{eqnarray*}
&&\mathop{\sum_{\lambda\vdash k}}_{\lambda=1^{m_1}2^{m_2}\cdots
}
\frac
{k!}{m_1!m_2!\cdots}\frac{1}{(2\pi\I)^{\ell(\lambda)}} \int_{(\alpha+\I\rr)^{\ell(\lambda)}}\,d\vec w \det
\biggl[\frac
{1}{w_a+\lambda_a-
w_b} \biggr]_{\lambda_1,\lambda_b=1}^{\ell(\lambda)}
\\
&&\qquad{}\times H(w_1,w_1+1,\ldots,w_1+
\lambda_1-1,\ldots,w_{\ell(\lambda
)},w_{\ell(\lambda)+1},
\ldots,w_{\ell(\lambda)+\lambda_{\ell
(\lambda)}-1})
\nonumber
\end{eqnarray*}
with
\[
H(z_1,\ldots,z_\ell)= \prod
_{a=1}^\ell e^{({t}/{2})z_a^2+z_ax}\sum
_{\sigma\in
S_\ell}\prod_a
\frac{1}{z_{\sigma(1)}+\cdots+z_{\sigma(a)}}\prod_{a>b}\frac{z_{\sigma(a)}-z_{\sigma(b)}-1}{z_{\sigma(a)}-z_{\sigma(b)}}.
\]
Using the same procedure as in \eqref{eq:mukns} to get rid of the
multinomial coefficient
$\frac{k!}{m_1!m_2!\cdots}$, the above turns into
\begin{eqnarray*}
&&\sum_{\ell=0}^k\frac{1}{\ell!}\mathop{
\sum_{m_1,\ldots,m_\ell
\geq1}}_{
m_1+\cdots+ m_\ell=k}\frac{1}{(2\pi\I)^\ell} \int
_{(\alpha+\I\rr)^{\ell}}\,d\vec w \det \biggl[\frac
{1}{w_a+m_a-w_b}
\biggr]_{\lambda_a,\lambda_b=1}^{\ell} \\
&&\qquad{}\times \prod_{a=1}^\ell
e^{({t}/{2})z_a^2+z_ax}
\\
&&\qquad{}\times H(w_1,w_1+1,\ldots,w_1+m_1-1,
\ldots,w_{\ell},w_{\ell
+1},\ldots ,w_{\ell}+m_\ell-1).
\end{eqnarray*}
In order to compute the sum over the symmetric group appearing in the
definition of $H$, we
will appeal to the following summation formula, which was used in \cite
{cal-led}.

%
\begin{lem}\label{lem:magic}\label{magic}
For $q_1,\ldots, q_N, \kappa\in\cc$,
\[
\sum_{\sigma\in S_N} \mu_{\vec{q}} (\sigma)\prod
_{a<b} \frac{
q_{\sigma
(a)} - q_{\sigma(b)}-i\kappa}{q_{\sigma(a)} - q_{\sigma(b)}} = \prod
_{a<b} \frac{ q_{a} +q_{b}+i\kappa}{q_{a} + q_{b}},
\]
where
\[
\mu_{\vec{q}} (\sigma):={q^{-1}_{\sigma(1)}
(q_{\sigma(1)} + q_{\sigma
(2)})^{-1} \cdots(q_{\sigma(1)} +
\cdots+ q_{\sigma(N) })^{-1} } {\prod_a
q_a}.
\]
\end{lem}

This identity was discovered and checked for small values of $N$ on
Mathematica by Le
Doussal and Calabrese. The formula can, in fact, be derived as a
suitable limit of an
analogous symmetrization identity proved in \cite{lee} in the context
of ASEP with flat
initial condition (see Lemma~2 in that paper).

Using the lemma, we obtain
\[
H(z_1,\ldots,z_\ell)=\prod
_{a=1}^\ell e^{
({t}/{2})z_a^2+z_ax}\prod
_{a<b}\frac{z_a+z_b-1}{z_a+z_b}.
\]
Replacing this formula above and doing some algebra leads directly to
\eqref{Eq:DBGHF-shifted-2} below, which as we will see next is another
way of
writing the solution given in Proposition~\ref{prop:deltaBoseFormulaU}
when all the
$x_a$'s are the same, thus proving \eqref{eq:deltaBose-lin}.

In what follows, we will turn our formula for the tilted half-flat
delta Bose gas (with all
$x_a$'s the same) into one in which all the integration contours
coincide. As we will see,
this alternative version of our half-flat formula is essentially
equivalent to the
formulas given in \cite{cal-led,leDoussalHF} [see \eqref
{Eq:DBGHF-shifted-2} and the
discussion that follows it]. We will argue afterwards (see Proposition~\ref{prop:DBGHF-shifted-2}), based on the convergence of ASEP to KPZ,
that this formula
does indeed give the half-flat SHE moments.

The first step is to deform the $z_a$ contours in \eqref{Eq:DBGHF} one
by one so that they
all coincide with $\alpha_k+\I\rr$. The arguments are similar to the
ones we used for ASEP
in Section~\ref{sec:mom}, so we only sketch them. We proceed similarly
to the proof of
Proposition~\ref{prop:muk}, now accounting for poles of the form
$z_a=z_b+1$ for $a>b$ and
computing the corresponding residues. Doing this in the case that all
$x_a$'s are
equal, using the symmetrization identity
\[
\sum_{\sigma\in S_k}\prod_{a>b}
\frac{z_{\sigma(a)}-z_{\sigma
(b)}-1}{z_{\sigma(a)}-z_{\sigma(b)}}=k!,
\]
which plays the role of \eqref{eq:symmmc} (and follows from suitably
rescaling it\footnote{It also corresponds to a
certain degeneration of the special case of the Hall--Littlewood
polynomial normalization
given in Section III.1 of \cite{macdonald}.}), and rewriting the sum
over partitions as in \eqref{eq:mukns}
yields the following formula for the moments of the delta Bose gas with
initial condition\vadjust{\goodbreak}
$v(0;x,\ldots,x)=e^{-\theta x}\mathbf{1}_{x\geq0}$ (here, and below,
$x$ is
repeated $k$ times in
the argument of $v$):
\begin{eqnarray*}
\label{Eq:DBGHF-shifted-pre} v(t;x,\ldots,x)&=&k!\sum_{\ell=0}^k
\frac{1}{\ell!}\mathop{\sum_{n_1,\ldots
,n_\ell,}}_{n_1+\cdots+n_\ell=k}
\frac{1}{(2\pi\I)^\ell} \int_{(\alpha+\I\rr)^\ell} \,d\vec{w} \det \biggl[
\frac{1}{w_a+n_a-w_b} \biggr]_{a,b=1}^\ell
\\
&&{}\times\bar H(w_1,\ldots,w_1+n_1-1,
\ldots,w_\ell,\ldots,w_\ell +n_\ell-1)
\end{eqnarray*}
with
\[
\bar H(z_1,\ldots,z_m)=\prod
_{a<b}\frac{z_a+z_b-1}{z_a+z_b}\prod_{a=1}^m
\frac{1}{z_a}e^{({t}/{2})
(z_a-\theta)^2+x(z_a-\theta)},
\]
where $\alpha>0$. Rewriting the result as in
the proof of Theorem~\ref{thmm:main-hf} yields (after some simplification)
%
%
\begin{eqnarray}\label{eq:vtxxx}
&&v(t;x,\ldots,x)\nonumber\\
&&\qquad=2^kk!\sum_{\ell=0}^k
\frac{1}{\ell!}\mathop{\sum_{n_1,\ldots,n_\ell,}}_{n_1+\cdots+n_\ell=k}
\frac{1}{(2\pi\I
)^\ell} \int_{(\alpha+\I\rr)^\ell} \,d\vec{w} \det \biggl[
\frac{1}{w_a+n_a-w_b} \biggr]_{a,b=1}^\ell
\nonumber\\
&&\qquad\quad{} \times\prod_a\frac{\Gamma(2w_a+n_a-1)}{\Gamma(2w_a+2n_a-1)}
\nonumber
\\[-8pt]
\\[-8pt]
\nonumber
&&\hspace*{26pt}\qquad\quad{}\times\exp\biggl\{\frac{1}2t \biggl[\frac{1}{3}n_a^3-\frac{1}{2}n_a^2+\frac
{1}{6}n_a+n_a(w_a-\theta
)^2+n_a(n_a-1)w_a\biggr]\nonumber\\
&&\hspace*{164pt}\qquad\quad{}+x \biggl[\frac{1}{2}n_a^2-\frac
{1}{2}n_a+n_a(w_a-\theta)
\biggr]\biggr\}\nonumber
\\
&&\qquad\quad{} \times\prod_{a<b}\frac{\Gamma(w_a+w_b+n_a-1)\Gamma(w_a+w_b+n_b-1)}{
\Gamma(w_a+w_b-1)\Gamma(w_a+w_b+n_a+n_b-1)}.\nonumber
\end{eqnarray}
Now we change variables $w_a\mapsto w_a-\frac{1}{2}(n_a-1)$ to obtain
%
%
\begin{eqnarray}
\label{Eq:DBGHF-shifted}&& v(t;x,\ldots,x)\nonumber\\
&&\qquad=2^kk!\sum
_{\ell=0}^k\frac{1}{\ell!}\mathop{\sum
_{n_1,\ldots,n_\ell\geq1}}_{n_1+\cdots+n_\ell=k}\frac{1}{(2\pi
\I )^\ell} \int
_{\alpha+({1}/{2})(n_1-1)+\I\rr} \,dw_1\cdots\\
&&\hspace*{95pt}\qquad\quad{}\times\int_{\alpha
+({1}/{2})(n_\ell-1)+\I\rr}
\,dw_\ell I_\theta(\vec w,\vec n)\nonumber
\end{eqnarray}
with
%
%
\begin{eqnarray}
\label{eq:Itheta}&& I_\theta(\vec w,\vec n)\nonumber\\
&&\qquad=
 \det \biggl[\frac{1}{w_a-w_b+({1}/{2})n_a+({1}/{2})n_b}
\biggr]_{a,b=1}^\ell\\
&&\qquad\quad{}\times \prod_a
\exp\biggl\{t\biggl [\frac{1}{24}n_a^3-\frac{1}{24}n_a+\frac
{1}{2}n_a(w_a-\theta)^2
\biggr]+xn_a(w_a-\theta)\biggr\}\nonumber
\\
&&\qquad\quad{}\times\prod_a\frac{\Gamma(2w_a)}{\Gamma(2w_a+n_a)}\nonumber\\
&&\qquad\quad{}\times \prod
_{a<b}\frac{\Gamma(w_a+w_b+({1}/{2})(n_a-n_b))\Gamma
(w_a+w_b-({1}/{2})(n_a-n_b))}{
\Gamma(w_a+w_b-({1}/{2})(n_a+n_b))\Gamma(w_a+w_b+({1}/{2})(n_a+n_b))}.
\nonumber
\end{eqnarray}
The last step is to shift back the $w_a$ contours from $\alpha+\frac
{1}{2}(n_a-1)+\I\rr$ to
$\alpha+\I\rr$. As we will see, we will not cross any poles as we do
this. To be more
precise, we begin by moving the $w_1$ contour from $\alpha+\frac
{n_1-1}2+\I\rr$ to
$\alpha+\I\rr$. There are three types of possible singularities, the
first from the Cauchy determinant and
the other two from the Gamma functions:
\begin{longlist}[(1)]
\item[(1)]$w_1=w_b-\frac{1}2 (n_1+n_b)$ for $b>1$.
\item[(2)]$w_1=-\ell$ for $\ell\in\zz_{\geq0}$.
\item[(3)]
$w_1=-w_b\pm\frac{1}2(n_1-n_b)-\ell$ for $\ell
\in
\zz_{\geq0}$ and $b>1$.
\end{longlist}
The first two types of singularity lie to the left of the origin,
whereas our deformation
region lies entirely to the right of the origin. Turning to (3), both
singularities may or may not lie inside the deformation region, but in
any case the
singularity is removable: the simple pole coming from the numerator
cancels with the zero of the denominator since $w_1=-w_b\pm\frac{1}2(n_1-n_b)-\ell$ implies
$w_1+w_b\mp\frac{1}2(n_1-n_b)=-\ell\in\zz_{<0}$, which is a zero of
$\frac{1}{\Gamma(\cdot)}$.

It remains to show that, having moved $w_1,\ldots,w_{j-1}$ from their
respective starting
points to $\alpha+\I\rr$, we do not incur any residues when moving
$w_j$ from
$\alpha+\frac{n_j-1}2+\I\rr$ to $\alpha+\I\rr$. The argument is
analogous to the case
$w_1$ and is left to the reader.

This leads to the following.

%
\begin{prop}\label{prop:DBGHF-shifted-2}
For the delta Bose gas with tilted half-flat initial condition
$v_0(t;\vec
x)=\prod_ae^{-\theta x_a}\mathbf{1}_{x_a\geq0}$ we have
%
%
\begin{equation}
\label{Eq:DBGHF-shifted-2} v(t;x,\ldots,x)=2^kk!\sum
_{\ell=0}^k\frac{1}{\ell!}\mathop{\sum
_{n_1,\ldots,n_\ell\ge1}}_{n_1+\cdots+n_\ell=k}\frac{1}{(2\pi\I
)^\ell} \int
_{(\alpha+\I\rr)^k}\,d\vec w\, I_\theta(\vec w;\vec n)
\end{equation}
with $I_\theta$ given by \eqref{eq:Itheta}. Moreover, in the pure
half-flat initial
condition corresponding to $\theta=0$, the same identity holds for the
moments of the
SHE, that is,
\[
\ee^{{\mathrm{ h\mbox{-}fl}}} \bigl[Z(t,x)^k \bigr]=2^kk!\sum
_{\ell
=0}^k\frac{1}{\ell !}\mathop{\sum
_{n_1,\ldots,n_\ell\ge1}}_{n_1+\cdots+n_\ell=k}\frac{1}{(2\pi\I
)^\ell} \int
_{(\alpha+\I\rr)^k}\,d\vec w\, I_0(\vec w;\vec n).
\]
\end{prop}

In \cite{cal-led}, the authors compute a formal series\footnote{As
written, the computation
in \cite{cal-led} is only formal, since in view of (88) in their paper
the series given
in their formula (40) is clearly divergent. Nevertheless, their
computation implicitly
leads to a formula like \eqref{Eq:DBGHF-shifted-2} in view of their
formula (39).} for
the generating function of $Z(t,x)$ using the explicit basis of
eigenfunctions of
the delta Bose gas \cite{liebLiniger,mcGuire}. The generating
function is expanded in
the ``number of strings'', which essentially corresponds to the
parameter $\ell$ in
\eqref{Eq:DBGHF-shifted-2} (the ``strings'' essentially correspond to
$n_1,\ldots,n_\ell$,
and index the eigenfunctions). The coefficients in this expansion are
given in their
formula (88), and one can check that, as expected, that formula
coincides essentially with
\eqref{Eq:DBGHF-shifted-2}. By this we mean that, for fixed
$n_1,\ldots
,n_\ell$, the
summand in \eqref{Eq:DBGHF-shifted-2} coincides\footnote{A diligent
reader will notice two
minor differences between \eqref{Eq:DBGHF-shifted-2} and the formula in
\cite{cal-led}. First, the formulas differ by a factor of $\prod_a(-1)^{n_a}$ in the
integrand, reflecting the fact that their generating function
computation is implicitly
calculating $\ee[(-Z)^k]$, as opposed to $\ee[Z^k]$ [see (40) in their
paper]. Second,
one needs to replace $t$ by $2t$ in our formula to recover theirs.
This is because in
their definition of the SHE the Laplacian term lacks the prefactor
$\frac{1}{2}$ [see (8) in
their paper and compare with \eqref{SHE}].} with the summand on the
right-hand side of (88) in~\cite{cal-led} with $n_s=\ell$ and $m_a=n_a$ for $a=1,\ldots,n_s$. This
correspondence is
consistent with (39) in their paper. See also \cite{leDoussalHF}.

\begin{pf*}{Proof of Proposition~\ref{prop:DBGHF-shifted-2}}
The delta Bose gas case follows directly from the above discussion.
The formula in the
case of the moments of the half-flat SHE can be recovered directly as
a weakly
asymmetric limit of the half-flat ASEP moment formula given in Theorem~\ref{thmm:main-hf}. Let us briefly sketch how this is done.

Recall from \eqref{eq:NflN0} that (for the half-flat case)
$h(t,x)=2N_x(t)-x$. According to
the WASEP scaling theory (see \cite{berGiaco}), if $\gamma=q-p=\ep
^{1/2}$ and we let
%
%
\begin{eqnarray}\label{eq:nnulambda}
\nnu_\ep&=&1-2\sqrt{pq}=\frac{1}2\ep+\frac{1}8
\ep^2+\mathcal{O}\bigl(\ep ^3\bigr),
\nonumber
\\[-8pt]
\\[-8pt]
\nonumber
\lambda_\ep&=&\frac{1}2\log\biggl(\frac{q}{p}\biggr)=
\ep^{1/2}+\frac{1}3\ep ^{3/2}+\mathcal{O}\bigl(
\ep^{5/2}\bigr),
\end{eqnarray}
then
\[
h_\ep(t,x):=\lambda_\ep h\bigl(\ep^{-3/2}t/
\gamma,\ep^{-1}x\bigr)-\nnu_\ep \ep^{-2}t
\]
converges to the Cole--Hopf solution $h(t,x)$ of KPZ starting with
$h(0,x) = 0$ for
$x>0$ and $h(0,x)=\infty$ for $x<0$, which in view of \eqref{eq:CH}
corresponds to $Z(0,x)=\mathbf{1}_{x\geq0}$.
Translating back to $N_x(t)$, and in view of \eqref{eq:nnulambda}, we have
\begin{eqnarray*}
N_{\ep^{-1}x}\bigl(\ep^{-2}t\bigr) &=& \frac{1}{2\lambda_\ep}
\bigl[h_\ep (t,x)+\nnu _\ep\ep^{-2}t \bigr]+
\frac{1}2\ep^{-1}x \\
&\approx&\frac{1}2\ep^{-1/2}h(t,x)-
\frac{1}{48}\ep^{-1/2}t+\frac{1}4\ep ^{-3/2}t+
\frac{1}2\ep^{-1}x.
\end{eqnarray*}
Therefore, since $\log(\tau)\approx-2\ep^{1/2}$, we deduce that
\[
\tau^{N_{\ep^{-1}x}(\ep^{-2}t)-
({1}/4)\ep^{-3/2}t-({1}/2)\ep^{-1}x}
\approx e^{-h(t,x)-({1}/{24})t}=e^{-({1}/{24})t}Z(t,x).
\]
This, together with tightness of the moments (which, e.g., can
be obtained by
adapting the arguments in Section~2.15 of \cite{quastelCDM}), gives
%
%
\begin{equation}\qquad
\label{eq:momASEPSHE} \ee^{{\mathrm{ h\mbox{-}fl}}}
\bigl[\tau^{k (N_{\ep^{-1}x}(\ep
^{-2}t)-({1}/4)\ep
^{-3/2}t-({1}/2)\ep^{-1}x+({1}/{24})t )} \bigr]
\mathop{\rightarrow}\limits_{\ep\to0}\ee^{{\mathrm{ h\mbox{-}fl}}} \bigl[Z(t,x)^k \bigr].
\end{equation}

Now in view of \eqref{eq:stpt1}, the left-hand side of \eqref
{eq:momASEPSHE} is given by
%
%
\begin{eqnarray}
\label{eq:momASEPSHEf} &&k_\tau!\sum_{\ell=0}^{k}
\frac{1}{\ell!}\mathop{\sum_{n_1,\ldots
,n_\ell
\geq1}}_{n_1+\cdots+n_\ell=k}
\frac{1}{(2\pi\I)^{\ell}}\int_{\gamma_{-1,0}^\ell}\,d\vec w \det \biggl[
\frac{-1}{w_a\tau^{n_a}-w_b} \biggr]_{a,b=1}^{\ell}
\nonumber\\
&&\qquad{}\times\prod_ae^{({1}/{24})n_at}
\tau^{-(({1}/4)\ep^{-3/2}t-({1}/2)\ep^{-1}x)n_a} \mff(w_a;n_a)\mfg(w_a;n_a)
\\
&&\qquad{}\times\prod_{a<b}\mfh(w_a,w_b;n_a,n_b).
\nonumber
\end{eqnarray}
Observe that the two sums are finite, so in order to obtain
\eqref{Eq:DBGHF-shifted-2} it is enough to show that the multiple
integral converges to
$I_0(\vec w;\vec n)$. This results in a relatively simple problem in asymptotic
analysis. The starting point for the critical point analysis is to
consider the product
$\tau^{-(({1}/4)\ep^{-3/2}t-({1}/2)\ep^{-1}x)n_a}\mff(w_a;n_a)$, which
is given
by
\begin{eqnarray*}
&&(1-\tau)^{n_a} \exp\biggl\{\biggl(\frac{1}{1+w_a}-\frac{1}{1+\tau^{n_a}w_a}-\frac{1}4\log(\tau)\biggr)
\ep^{-3/2}t\\
&&\hspace*{42pt}\qquad{}-
\frac{1}2\ep^{-1}x\log(\tau)n_a+\log\biggl(\frac{1+\tau
^{n_a}w_a}{1+w_a}\biggr)\bigl(\ep^{-1}x-1\bigr)\biggr\}.
\end{eqnarray*}
Scaling $w_a$ near 1 through the change of variables $w_a\mapsto
1-(1-\tau)\tilde w_a$,
the exponent above can be written as $(\frac{1}6n_a^3+\frac{1}2n_a^2\tilde
w_a+\frac{1}2n_a\tilde w_a^2)t+(\frac{1}2n_a^2+n_a\tilde w_a)x+\mathcal
{O}(\ep^{1/2})$. The
change of variables (for all the $\ell$ variables) gives a prefactor of
$(-1)^\ell(1-\tau)^\ell$, while the factor $\prod_a(1-\tau)^{n_a}$
coming from the above
product turns into $(1-\tau)^k$. We leave it to the reader to verify
that, with this
scaling, $\det [\frac{-1}{w_a\tau^{n_a}-w_b}
]_{a,b=1}^{\ell}
\approx(1-\tau)^{-\ell}\det [\frac{1}{\tilde w_a+n_a-\tilde
w_b} ]_{a,b=1}^{\ell}$,
$\prod_a\mfg(w_a;n_a)\approx2^k(1-\tau)^{-k}\prod_a\frac{\Gamma
(2\tilde w_a+n_a)}{
\Gamma(2\tilde w_a+2n_a)}$ and
\[
\prod_{a<b}\mfh(w_a,w_b;n_a,n_b)\approx\prod_{a<b}\frac{\Gamma
(\tilde
w_a+\tilde
w_b+n_a)\Gamma(\tilde w_a+\tilde w_b+n_b)}{\Gamma(\tilde w_a+\tilde
w_b+n_a+n_b)
\Gamma(\tilde w_a+\tilde w_b)}.
\]
Note that near the critical point
$w_a=1$ the contour
$\gamma_{-1,0}$ turns into $\I\rr$, negatively oriented.\vadjust{\goodbreak} Introducing
an additional
factor $(-1)^\ell$ to flip the orientation of the resulting contour,
we deduce from the
above estimates that \eqref{eq:momASEPSHEf} is approximately
\begin{eqnarray*}
&&2^kk!\sum_{\ell=0}^k
\frac{1}{\ell!}\mathop{\sum_{n_1,\ldots
,n_\ell,}}_{
n_1+\cdots+n_\ell=k}
\frac{1}{(2\pi\I)^\ell} \int_{(\I\rr)^\ell} \,d\vec{\tilde w} \det \biggl[
\frac{1}{\tilde
w_a+n_a-\tilde w_b} \biggr]_{a,b=1}^\ell\\
&&\qquad{}\times\prod
_a\frac{\Gamma
(2\tilde
w_a+n_a)}{\Gamma
(2\tilde w_a+2n_a)}
\\
&&\qquad{}\times\prod_a\exp\biggl\{t \biggl[\frac{1}6n_a^3+\frac{1}2n_a^2\tilde
w_a+\frac{1}2n_a\tilde w_a^2+\frac{1}{24}n_a \biggr]+x \biggl[\frac{1}2n_a^2+n_a
\tilde w_a \biggr]\biggr\} \\
&&\qquad{}\times\prod
_{a<b}\frac{\Gamma(\tilde w_a+\tilde w_b+n_a)\Gamma(\tilde
w_a+\tilde w_b+n_b)}{\Gamma(\tilde w_a+\tilde w_b)\Gamma(\tilde
w_a+\tilde w_b+n_a+n_b)}.
\nonumber
\end{eqnarray*}
Turning this into a rigorous proof involves estimating the integrand
away from the
critical point in order to show that the only contribution from the
integral that
survives in the limit is that near $w_a=1$. This is not hard to do in
this case because
we do not need an estimate which is uniform in $\ell$ (which is the
basic source of difficulty in
turning the calculations of the \hyperref[sec:app-2to1]{Appendix} into a
rigorous proof), so we
will leave the details to the reader. Now changing variables $\tilde
w_a\mapsto\bar
w_a-\frac{1}2$ turns the above formula into \eqref{eq:vtxxx} (with
$\alpha=\frac{1}2$), which by
\eqref{Eq:DBGHF-shifted} gives the desired result.
\end{pf*}

\begin{appendix}\label{sec:app-2to1}
\section*{Appendix: Asymptotics for half-flat ASEP and the \texorpdfstring{Airy$_{2\to1}$}{Airy2to1}
marginals}

In this section, we provide a formal critical point analysis of the
long-time asymptotics
of the $\tau$-Laplace transform of $\tau^{N_{x}(t)}$ in the half-flat
case which, in
view of \eqref{eq:eepp-hf}, gives the asymptotic distribution of the
fluctuations of the
height function $h(t,x)$.

More precisely, our derivation will provide a nonrigorous confirmation
of the
conjectured asymptotics
%
%
\setcounter{equation}{0}
\begin{eqnarray}
\label{eq:2to1-conj}&& \lim_{t\to\infty}\pp^{\mathrm{ h\mbox{-}fl}} \biggl(
\frac
{h(t/(q-p),t^{2/3}x)-({1}/2)t-t^{1/3}x^2\mathbf{1}_{x\leq0}}{t^{1/3}}\geq-r
 \biggr)
 \nonumber
 \\[-8pt]
 \\[-8pt]
 \nonumber
 &&\qquad=\pp\bigl({\mathcal{A}}_{2\to1}
\bigl(2^{-1/3}x\bigr)\leq2^{1/3}r\bigr),
\end{eqnarray}
where ${\mathcal{A}}_{2\to1}$ is the Airy$_{2\to1}$ process. For
background on this
process and more details
about this conjecture, see \cite{quastelRem-review}.

Our starting point is the formula for the $e_\tau$-Laplace transform\break  of
$\tau^{N_{x}(t)}$
given in Theorem~\ref{thmm:gen-hf}, where we take
$\zeta=\break -\tau^{-({1}/{4})t-({1}/{2})t^{2/3}x+
({1}/{2})t^{1/3}r-({1}/{4})t^{1/3}x^2\mathbf{1}_{x\leq0}}$ and let
$\tilde r=r-\frac{1}{2}x^2\mathbf{1}_{x\leq0}$:
%
%
\begin{eqnarray}
\label{eq:gen-hf-app} &&\ee^{\mathrm{ h\mbox{-}fl}} \bigl[e_\tau\bigl(-
\tau^{N_{t^{2/3}x}(t/\gamma
)-({1}/{4})t-({1}/{2})t^{2/3}x+({1}/{2})t^{1/3}\tilde
r}\bigr) \bigr]\nonumber
\\
&&\qquad=\sum_{k=0}^{\infty}\frac{1}{k!}
\frac{1}{(2\pi\I)^{2k}}\int_{(\delta+\I\rr
)^k}\,d\vec s \int_{\gamma_{-1,0}^k}\,d
\vec w \det \biggl[\frac{-1}{w_a\tau
^{s_a}-w_b} \biggr]_{a,b=1}^{k}
\nonumber
\\[-8pt]
\\[-8pt]
\nonumber
&&\qquad\quad{}\times\prod_a\tau^{-[({1}/{4})t+
({1}/{2})t^{2/3}x-
({1}/{2})t^{1/3}r]s_a}\tilde\mff
(w_a;s_a)\mfg(w_a;s_a)\\
&&\qquad\quad{}\times\prod
_{a<b}\mfh(w_a,w_b;s_a,s_b)
\nonumber
\end{eqnarray}
for $\delta\in(0,1)$, $\mff$, $\mfg$ and $\mfh$ as in \eqref
{eq:germans-hf}, and
with $\tilde\mff$ defined as $\mff$ with $t$ replaced by $t/\gamma$
(recall that $\gamma=q-p$).

We will perform a formal critical point analysis on the right-hand
side. The reason the
limit is not rigorous is that so far we have not been able to control
the double product
$\prod_{a<b}\mfh(w_a,w_b;s_a,s_b)$ on the part of the contour away
from the critical
point, nor find an alternative contour where this can be done.\footnote
{More precisely,
$|\mfh(w_a,w_b;s_a,s_b)|$ can be bounded uniformly by some constant
$C$, but this
constant is necessarily larger than one. This yields an estimate of
the form $C^{k^2}$
for some $C>1$, which is too big (note that $\sum_{k\geq0}\frac
{1}{k!}C^{k^2}$ is
divergent). Therefore the rigorous asymptotics remains an interesting open
problem. Observe that if the double product could be turned into a
determinant [as
happens for the first double product in \eqref{Eq:NestedPlusSmall},
which turns into the
determinant in \eqref{eq:gen-hf-app}], then this problem would
disappear, because by
Hadamard's bound our estimate on $|\mfh(w_a,w_b;s_a,s_b)|$ would
essentially yield a
factor $C^kk^{k/2}$, which is small enough for our purposes.} The
derivation here is
done to clarify the algebraic structure of the expansion around the
critical point where
one sees the Airy crossover distributions.

The leading order (in $t$) factor in the integrand comes $\tilde\mff
(w_a,s_a)$ and the factor
$\tau^{-({1}/{4})t}$, and can be written as
$\prod_a\exp [t  (\frac{1}{1+w_a}-\frac{1}{1+\tau
^{s_a}w_a}-\frac{1}4s_a\log(\tau) ) ]$.
One can verify that the only critical point of
$\frac{1}{1+w}-\frac{\tau}{1+\tau^{s}w}-\frac{1}4s\log(\tau)$ occurs at
$(w,s)=(1,0)$. Moreover, the Hessian of this function vanishes at this
point, while the
third-order partial derivatives are not all 0, which suggests a
$t^{1/3}$ scaling. On the
other hand, this suggests that the $w_a$ contour should be chosen to
cross the line
$\rr_{\geq0}$ at $w_a=1$.
In view of this we change variables as follows:
%
%
\begin{equation}
\label{scaling123} w_a=1+t^{-1/3}\tw_a,\qquad
s_a=-\frac{1}{\log(\tau)}t^{-1/3}\ts_a.
\end{equation}
We will need the following lemma.

%
\begin{lemm}
\label{Lem:PochAsymptotics}
Let $a\in\rr$, $\ell\in\zz$ and $k\in\zz_{>0}$. If $-\ell=kj$
for some
$j\in\zz_{\geq0}$ (i.e., $\ell=0$ or $k$ is a factor of $-\ell$) then,
as $\varepsilon\to0$,
\[
\label{Eq:PochAsymptotics} \bigl(\tau^{\ell}(1+\varepsilon a);\tau^k
\bigr)_\infty = -\varepsilon a\mathop{\prod_{n=0}}_{n\neq j}^\infty
\bigl(1-\tau^{kn+\ell
} \bigr) + \mathcal{O}\bigl(\ep^2\bigr).
\]
On the other hand, if $\ell\neq0$ and $k$ is not a factor of $-\ell$
then, as $\varepsilon\to0$,
\[
\label{Eq:PochAsymptotics2} \bigl(\tau^\ell(1+\varepsilon a);\tau^k
\bigr)_\infty = \bigl(\tau^\ell;\tau^k
\bigr)_\infty \Biggl[1-\varepsilon a\sum_{m=0}^\infty
\frac{\tau^{km+\ell}}{1-\tau^{km+\ell}} \Biggr] + \mathcal{O}\bigl(\ep^2\bigr).
\]
\end{lemm}

\begin{pf}
In the first case, we have $\tau^{kj+\ell}=1$ so
\begin{eqnarray*}
\bigl(\tau^{\ell}(1+\varepsilon a);\tau^k \bigr)_\infty
&= &\prod_{n=0}^\infty \bigl[1- (1+\varepsilon
a)\tau^{kn+\ell} \bigr] \\
&= &\bigl[1- (1+\varepsilon a) \bigr] \mathop{\prod
_{n=0}}_{n\neq
j}^\infty \bigl[1- (1+
\varepsilon a)\tau^{kn+\ell} \bigr]
\\
&=& -\varepsilon a \mathop{\prod_{n=0}}_{n\neq j}^\infty
\bigl(1-\tau ^{kn+\ell} \bigr) + \mathcal{O}\bigl(\ep^2\bigr).
\end{eqnarray*}
In the second case, we have
\begin{eqnarray*}
\bigl(\tau^\ell(1+\varepsilon a);\tau^k
\bigr)_\infty &=& \prod_{n=0}^\infty
\bigl[\bigl(1- \tau^{kn+\ell}\bigr)-\varepsilon a\tau^{kn+\ell} \bigr]
\\
&=&\prod_{n=0}^\infty \bigl[1-
\tau^{kn+\ell} \bigr]-\varepsilon a\sum_{m=0}^\infty
\tau^{km+\ell}\mathop{\prod_{n=0}}_{n\neq
m}^\infty
\bigl[1-\tau^{kn+\ell} \bigr]+\mathcal O\bigl(\ep^2\bigr)
\\
&\approx& \bigl(\tau^\ell;\tau^k \bigr)_\infty-
\varepsilon a \bigl(\tau^\ell;\tau^k \bigr)_\infty
\sum_{m=0}^\infty\frac
{\tau^{km+\ell
}}{1-\tau^{km+\ell}}.\quad\qed
\end{eqnarray*}
\noqed\end{pf}

The scaling \eqref{scaling123} leads to the following asymptotics:
\begin{eqnarray*}
&&\frac{\pi}{\sin(-\pi s_a)}\frac{1+w_a}{1+\tau^{s_a}w_a}(1-\tau )^{s_a}
\approx
\frac{\log(\tau)t^{1/3}}{\ts_a}, \\
&&\frac{1}{w_a\tau^{s_a}-w_b} \approx \frac{\tau^{-1}t^{1/3}}{\tw_a-\tw_b-\ts_a},
\\
&&t \biggl[\frac{1}{1+w_a}-\frac{1}{1+\tau^{s_a}w_a}-\frac{1}4s_a
\log (\tau )+\frac{1}2t^{1/3}s_ar\log(\tau) \biggr]
\\
&&\qquad\approx\frac{1}{48}\bigl(\ts_a^3-3
\ts_a^2\tw_a+3\ts_a
\tw_a^2\bigr)-\tilde r\ts _a,
\\
&&\biggl(\frac{1+\tau^{s_a}w_a}{1+w_a}
\tau^{-({1}/{2})s_a} \biggr)^{t^{2/3}x} \approx
\biggl(1-\frac{\ts_a^2-2\ts_a\tw_a}{8t^{2/3}} \biggr)^{t^{2/3}x}
\approx e^{-({1}/{8})(\ts_a^2-2\ts_a\tw_a)x},
\end{eqnarray*}
while, using Lemma~\ref{Lem:PochAsymptotics},
\begin{eqnarray*}
\frac{ (-w_a;\tau )_\infty}{ (-\tau^{s_a}w_a;\tau
 )_\infty}
\frac
{ (\tau^{2s_a}w_a^2;\tau )_\infty}{ (\tau
^{s_a}w_a^2;\tau )_\infty} &\approx&\frac{ (1+2(\tw_a-\ts_a)t^{-1/3};\tau )_\infty
}{ (1+(2\tw _a-\ts_a)t^{-1/3};\tau )_\infty}\\
&\approx&\frac{2(\ts_a-\tw_a)}{\ts_a-2\tw
_a},
\end{eqnarray*}
and similarly
\[ \frac{ (w_aw_b;\tau )_\infty (\tau
^{s_a+s_b}w_aw_b;\tau )_\infty}{
 (\tau^{s_a}w_aw_b;\tau )_\infty (\tau
^{s_b}w_aw_b;\tau )_\infty} \approx\frac{(\tw_a+\tw_b)(\tw_a+\tw_b-\ts_a-\ts_b)}{(\tw
_a+\tw_b-\ts
_a)(\tw_a+\tw_b-\ts_b)}.
\]
Additionally, there is a factor of $(-1)^kt^{-2k/3}(\tau/\log(\tau))^k$
coming from the
change of variables which, except for the $(-1)^k$, cancels exactly
with factors coming
out from the first line of the above list of asymptotics. To write the
limit choose first
$\delta=-t^{-1/3}/(2\log(\tau))$ in \eqref{eq:gen-hf-app} and deform
the $s_a$ contour so
that it departs the real axis at angles $\pm\pi/3$, and likewise deform
the $w_a$ contours
so that they go through $1$ and depart from that point at angles $\pm
\pi
/3$. The
limiting contours then become $\frac{1}{2}+\langle$ for $\ts_a$ and
$\langle$
for $\tw_a$, where
$\langle$ consists on two infinite rays departing 0 at angles $\pm\pi
/3$ (oriented with
increasing imaginary part) and thus using the above asymptotics in
\eqref{eq:gen-hf-app} we
obtain that the formal limit as $t\to\infty$ of $\ee [e_\tau
(-\tau
^{N_{t^{2/3}x}(t/(q-p))-({1}/{4})t-({1}/{2})t^{2/3}x+
({1}/{2})t^{1/3}\tilde
r}) ]$ is given by
\begin{eqnarray*}
F_x(\tilde r)&=&\sum_{k=0}^{\infty}
\frac{1}{k!}\frac{1}{(2\pi\I
)^{k}}\int_{({1}/{2}+\langle)^k}\,d\vec{\wt s}
\frac{1}{(2\pi\I)^{k }}\int_{(\langle)^k}\,d\vec{\wt w} \det \biggl[
\frac
{1}{\tw_b-\tw_a+\ts_a} \biggr]_{a,b=1}^{k}
\\
&&{} \times\prod_a\exp\biggl\{\frac{1}{48}(\ts_a^3-3\ts_a^2\tw_a+3\ts_a\tw
_a^2)-\frac{1}{2}\tilde r\ts_a-\frac{1}{8}(\ts_a^2-2\ts_a\tw_a)x\biggr\}\\
&&\hspace*{25pt}{}\times \frac{2(\ts_a-\tw_a)}{\ts_a(\ts_a-2\tw_a)}\\
&&{}\times \prod_{a<b}\frac{(\tw
_a+\tw
_b)(\tw_a+\tw_b-\ts_a-\ts_b)}{(\tw_a+\tw_b-\ts_a)(\tw_a+\tw
_b-\ts_b)}.
\end{eqnarray*}
Now we introduce the change of variables $\tw_a=u_a$ and $\ts
_a=u_a-v_a$. The $u_a$
contour is $\langle$, but we may freely deform it (thanks to the cubic
terms in the
exponent) to $1+\langle$. A priori $v_a$ depends on $u_a$, but again
one can check that
it can be deformed to $\rangle$, which is defined in the same way as
$\langle$ but departing the origin at
angles $\pm2\pi/3$. We obtain
\begin{eqnarray*}
\label{eq:MB-hf} F_x(\tilde r)&=&\sum_{k=0}^{\infty}
\frac{1}{k!}\frac{1}{(2\pi\I
)^{k}}\int_{(1+\langle)^k}\,d\vec{u}
\frac{1}{(2\pi\I)^{k}} \int_{(\rangle)^k}\,d\vec{v} \det \biggl[
\frac
{1}{u_b-v_a} \biggr]_{a,b=1}^{k}
\\
&&{}\times\prod_a\exp\biggl\{\frac{1}{48}(u_a^3-v_a^3)+\frac
{1}{8}(u_a^2-v_a^2)x-\frac{1}{2}(u_a-v_a)\tilde r\biggr\}
\frac{2v_a}{u^2_a-v^2_a}\\
&&{}\times\prod_{a<b}\frac{(u_a+u_b)(v_a+v_b)}{
(u_a+v_b)(v_a+u_b)}.
\nonumber
\end{eqnarray*}
Now we note that the determinant and the cross-product above simplify
into a single
determinant: using the Cauchy determinant formula
%
%
\[
\label{eq:cauchydet} \det \biggl[\frac{1}{x_a-y_b} \biggr]_{a,b=1}^k=
\frac{\prod_{a<b}(x_a-x_b)(y_b-y_a)}{\prod_{a,b}(x_a-y_b)},
\]
we have
\begin{eqnarray*}
&&\det \biggl[\frac{1}{u_b-v_a} \biggr]_{a,b=1}^k \prod
_{a<b} \frac{(
u_a+ u_b)( v_a+ v_b)}{( u_a+ v_b)( v_a+ u_b)} \\
&&\qquad= \frac{1}{\prod_a ( u_a-v_a)}
\prod_{a<b} \frac{( u_b^2- u_a^2)( v_a^2-
v_b^2)}{( u_a^2 - v_b^2)( u_b^2 - v_a^2)}
\\
&&\qquad= \frac{\prod_a ( v_a^2 - u_a^2)}{\prod_a ( v_a- u_a)} \det \biggl[ \frac{1}{u_a^2- v_b^2} \biggr]_{a,b=1}^k
= \prod_a(u_a+v_a)\det
\biggl[ \frac{1}{ u_b^2- v_a^2} \biggr]_{a,b=1}^k.
\end{eqnarray*}
Using this above, we get
\begin{eqnarray*}
F_x(\tilde r)& =& \sum_{k=0}^{\infty}
\frac{1}{k!}\frac{1}{(2\pi\I
)^{k}}\int_{(1+\langle)^k}\,d\vec{u}
\frac{1}{(2\pi\I)^{k }}\int_{(\rangle)^k}\,d\vec{v} \det \biggl[
\frac
{1}{u_b^2-v_a^2} \biggr]_{a,b=1}^{k}
\\
&&{}\times\prod_a\frac{2v_a}{u_a-v_a}e^{({1}/{48})(u_a^3-v_a^3)-
({1}/{8})(u_a^2-v_a^2)x-
({1}/{2})(u_a-v_a)\tilde
r}
\\
&=&\sum_{k=0}^{\infty}\frac{1}{k!}
\frac{1}{(2\pi\I)^{k}}\int_{(1+\langle)^k}\,d\vec{u}\\
&&{}\times \det \biggl[
\frac{1}{2\pi\I
}\int_{\rangle}\,dv \frac{2v}{u_a-
v}
\frac{e^{({1}/{48})u_a^3-({1}/{8})u_a^2x-({1}/{2})u_a\tilde
r}}{e^{({1}/{48})v^3-({1}/{8})v^2x-({1}/{2})v\tilde
r}} \frac{1}{u_b^2-v_a^2} \biggr]_{a,b=1}^k.
\end{eqnarray*}
This last expression is just the series expansion of a Fredholm determinant:
\[
F_x(\tilde r)=\det(I-K)_{L^2(1+\langle)}
\]
with
\begin{eqnarray*}
K\bigl(u,u'\bigr)&=&\frac{1}{2\pi\I}\int_{\rangle}\,dv
\frac{2v}{u-
v}\frac{e^{({1}/{48})u^3+({1}/{8})u^2x-({1}/{2})u_a\tilde
r}}{e^{({1}/{48})v^3-({1}/{8})v^2x-({1}/{2})v\tilde r}} \frac{1}{u'^2-v^2}
\\
&=&\frac{1}{2\pi\I}\int_0^\infty \,d\lambda\int
_{\rangle}\,dv \frac
{2v}{u'^2-v^2} \frac{e^{({1}/{48})u^3+({1}/{8})u^2x-u(\lambda+
({1}/{2})\tilde r)}}{e^{({1}/{48})v^3+({1}/{8})v^2x-v(\lambda+
({1}/{2})\tilde
r)}}.
\end{eqnarray*}
For more details on Fredholm determinants, see Section~2 of \cite
{quastelRem-review}.\vspace*{1pt}
Using the cyclic property of the Fredholm determinant, we
deduce that $F_x(\tilde r)=\det(I-\wt K)_{L^2([0,\infty))}$ with $\wt
K(\lambda,\lambda')=\frac{1}{(2\pi\I)^{2}}\int \,du \int \,dv \frac
{2v}{u^2-v^2}
\frac{e^{({1}/{48})u^3+({1}/{8})u^2x-u(\lambda+
({1}/{2})\tilde r)}}{e^{({1}/{48})v^3+({1}/{8})v_a^2x-v(\lambda
'+({1}/{2})\tilde
r)}}$\break 
and the same $u$ and $v$ contours. Scaling $u$ and $v$ by $2^{4/3}$ and
changing variables $\lambda\mapsto2^{-4/3}\lambda-\frac{1}{2}r$
and $\lambda'\mapsto2^{-4/3}\lambda'-\frac{1}{2}r$ finally yields
%
%
\begin{eqnarray}
\label{eq:limAiry2to1} &&\lim_{t\to\infty}\ee \bigl[e_\tau\bigl(-
\tau^{N_{t^{2/3}x}(t/\gamma
)(-{1}/{4})t-({1}/{2})t^{2/3}x+({1}/{2})t^{1/3}\tilde
r}\bigr) \bigr]
\nonumber
\\[-8pt]
\\[-8pt]
\nonumber
&&\qquad=\det \bigl(I-K^{2\to1}
\bigr)_{L^2([2^{1/3}r,\infty))}
\end{eqnarray}
with
\[
K^{2\to1}\bigl(\lambda,\lambda'\bigr)=\frac{1}{(2\pi\I)^{2}}
\int_{1+\langle}\,du \int_{\rangle
}\,dv
\frac{2v}{u^2-v^2} \frac{e^{({1}/{3})u^3+2^{-1/3}u^2x-u(\lambda-2^{-2/3}x^2\mathbf
{1}_{x\leq 0})}}{e^{({1}/{3})v^3+2^{-1/3}v^2x-v(\lambda
'-2^{-2/3}x^2\mathbf{1}_{x\leq0})}}.
\]
The $u$ and $v$ contours can be easily deformed to match those
appearing in the kernel
inside the Fredholm determinant which gives the finite dimensional
distributions of the
Airy$_{2\to1}$ process, see \cite{bfs}. Comparing with that formula, we
deduce that the
right-hand side of~\eqref{eq:limAiry2to1} equals $\pp({\mathcal
{A}}_{2\to1}
(2^{-1/3}x)\leq
2^{1/3}r)$ which,
in view of~\eqref{eq:eepp-hf}, finishes our formal derivation of
\eqref
{eq:2to1-conj}.
\end{appendix}
\section*{Acknowledgments}
The authors would like to thank Alexei Borodin, Ivan Corwin and Pierre
Le Doussal for
numerous discussions about the results in this paper. They would also
like to thank the
referees for their thorough review and comments.

%





\printaddresses
\end{document}